\documentclass[12pt,leqno]{article}
\usepackage{amsmath,amssymb,amsthm,amscd}
\newcommand{\doublespace}{
   \renewcommand{\baselinestretch}{1.2}
   \large\normalsize}

\doublespace

\let\ftilde\tilde %Plain
\renewcommand{\tilde}[1]{\mathbf{\ftilde{\mathnormal{#1}}}}

\newcommand{\C}{\mathbb{C}}
\newcommand{\Z}{\mathbb{Z}}

\newcommand{\N}{\mathbb{N}}
\newcommand{\Q}{\mathbb{Q}}

\newcommand{\unit}{\mathbf{1}}

\newcommand{\End}{\operatorname{End}}
\newcommand{\Res}[1]{\underset{#1=0}{\operatorname{Res}}\,}

\newcommand{\frakh}{\mathfrak{h}}

\newcommand{\id}{\operatorname{id}}

\newcommand{\calh}{\mathcal{H}}
\newcommand{\cali}{\mathcal{I}}
\newcommand{\calm}{\mathcal{M}}

\newcommand{\calw}{\mathcal{W}}
\newcommand{\cals}{\mathcal{S}}

\newcommand{\wt}[1]{\operatorname{wt}\,(#1)}
\newcommand{\NO}{\,{\raise0.25em\hbox{$\mathop{\hphantom{\cdot}}%
\limits^{_{\circ}}_{^{\circ}}$}}\,}

\def \End{{\rm End}}

\def \mod{{\rm mod}}
\def \<{\langle} 
\def \>{\rangle} 

\def \a{\alpha }

\def \l{\lambda }
\catcode`\@=11

\@addtoreset{equation}{subsection}
\makeatother
\catcode`\@=\active

\theoremstyle{plain}
 \newtheorem{theorem}{Theorem}[subsection]
	\newtheorem{corollary}[theorem]{Corollary}
	\newtheorem{lemma}[theorem]{Lemma}
	\newtheorem{proposition}[theorem]{Proposition}
	
\theoremstyle{remark}
 
 \newtheorem{remark}[theorem]{Remark}

\begin{document}

%title
\begin{large}
\begin{center}
Classification of irreducible modules\\
for the vertex operator algebra 
$M(1)^+$\, II: higher rank
\end{center}
\end{large}

\vskip 10ex

\begin{center}
Chongying Dong\footnote{Supported by NSF grant 
DMS-9700923 and a research grant from the Committee on Research, UC Santa
Cruz.}\\
%dong@math.ucsc.edu\\
Department of Mathematics, University of California, 
Santa Cruz, CA 95064, U.S.A\\

\vskip 2ex
Kiyokazu Nagatomo\footnote{Supported in part by  Grant-in-Aid for
Scientific Research,
the Ministry of Education, Science and Culture.}\\
%nagatomo@math.sci.osaka-u.ac.jp\\
Department of Mathematics, Graduate School of Science,
Osaka University\\
Osaka, Toyonaka 560-0043, Japan
\end{center}

\vskip 10ex

%Abstract

\begin{small}

\noindent
\textbf{Abstract}:
The vertex operator algebra $M(1)^+$ is the fixed point set of
free bosonic vertex operator algebra $M(1)$ of rank $\ell$ under the 
$-1$ automorphism. All irreducible modules for $M(1)^+$ 
are classified in this paper for any $\ell.$ 
\end{small}

\section{Introduction}\label{introduction}

This is the third paper in studying $\theta$-orbifold models associated to
lattice vertex operator algebras $V_L$ for even integral lattices $L$
where $\theta$ is an automorphism of $V_L$  of order 2 lifted from the
$-1$ isometry of $L.$ The $V_L$ contains the rank $\ell$ free bosonic
vertex  operator algebra $M(1)$ and the automorphism $\theta$ preserves
$M(1).$   In \cite{DN1} we studied  the orbifold model $M(1)^+$ which is
the $\theta$-fixed point set of $M(1)$ and classified all the inequivalent
irreducible modules by determining associated Zhu's algebra
$A(M(1)^+)$ explicitly  in the case of rank one.
The results and the method developed in
\cite{DN1} were effectively used in \cite{DN2} to get the classification
result for the inequivalent irreducible modules for the charge conjugation
orbifold model, which is the $\theta$-invariants of a
lattice vertex operator algebra $V_L$ for a rank one lattice $L.$
In this paper, we investigate the $\theta$-orbifold model $M(1)^+$ for 
arbitrary rank $\ell$ free bosonic vertex operator algebra $M(1)$ and classify
the irreducible modules for $M(1)^+.$ The results in this paper are
expected to be used to study the representation theory for a vertex
operator algebra $V_L^+$ which is the $\theta$-invariants of $V_L$ for
a lattice $L$ of rank $\ell.$

The free bosonic vertex operator algebra ${\cal H}=M(1)$ of 
rank $\ell$ (cf.~[FLM]) 
is an affine vertex operator algebra associated to an $\ell$-dimensional 
abelian Lie algebra $\frakh$ (see Subsection 2.2 below). The map
$\theta:\frakh\longrightarrow \frakh$ defined by $\theta(h)= -h$ induces
a vertex operator algebra automorphism denoted by the same symbol $\theta$.
Then the fixed point set ${\cal H}^+$  of $\theta$ is a simple 
vertex operator subalgebra of $\calh$. It is well known that
all the irreducible modules for $\calh$ are exhausted by Fock representation
$M(1,\lambda)$ for the affine algebra $\hat \frakh$ with the highest weight 
$\lambda\in \frakh.$ 
As a module for $\calh^+$,  $M(1,\lambda)$ and $M(1,-\lambda)$ are
isomorphic and irreducible if $\lambda \neq 0.$
But $M(1,0)={\cal H}$ decomposes into its irreducible
components $\calh = \calh^+\oplus\calh^-$ where $\calh^{\pm}$ are
the eigenspaces of $\theta.$ 
 One of the features of orbifold
models is the existence of extra irreducible modules which
come from the twisted sectors. The $\calh$ has 
exactly one irreducible $\theta$-twisted module
$\calh(\theta)$ with the $\theta$ action, which gives rise to two 
inequivalent irreducible modules
$\calh(\theta)^{\pm}$ for $\calh^+$ where $\calh(\theta)^{\pm}$ are
the eigenspaces of $\theta.$ The main result in this paper is that
$M(1,\lambda)$ ($\lambda\ne 0$), $\calh^{\pm}$ and $\calh(\theta)^{\pm}$ 
are all inequivalent irreducible $\calh^{+}$-modules.

In \cite{Z}, Zhu introduced an associative algebra $A(V)$ for any
vertex operator algebra $V,$ which gives a lot of
information on $V$ as far as the representation theory concerns. For
instance, there is a one to one correspondence between the set of
equivalence classes of irreducible modules for the associative algebra
$A(V)$ and the set of equivalence classes of irreducible admissible 
modules for $V$. This fact
has been used to classify the irreducible modules for
affine vertex operator algebras \cite{FZ}, Virasoro vertex operator algebras
\cite{W}, lattice vertex operator algebras \cite{DLM3},
$M(1)^+$ in the case $\ell=1$ [DN1] and the $-1$ orbifold vertex
operator algebra associated to the rank 1 lattice [DN2]. 
This idea was
developed further in [DLM2] to deal with twisted representations 
and the $\theta$-twisted modules for lattice vertex operator algebras $V_L$
were classified along this line \cite{DN3} .

The classification result in this paper is also achieved by using Zhu's 
algebra. The strategy is to determine Zhu's algebra $A(\calh^+)$
and to find a set of good generators and their relations.
The determination of Zhu's algebra is not only related to
the representation theory but also 
the structure theory for a given vertex operator algebra. 
For example we found a Poincar\'{e}-Birkhoff-Witt type theorem
for a $\calh^+$ in \cite{DN1} in the case $\ell=1.$  
So investigation of Zhu's algebra sheds light on the hidden structure of VOA's.

It is worth pointing out that there is a  main difference between rank one case
and the others.  Zhu's algebra for rank one case is
commutative but is not for higher rank case. For instance the top level
of the module $\calh^-$ is $\ell$-dimensional. The algebra structure
of $A(\calh^+)$ and ideas given in \cite{DN1} in the rank one case 
are very helpful but not enough to attack the higher
rank case. To overcome  the difficulty arising from noncommutativity
of $A(\calh^+)$ we introduce an ideal $\cali$ which is isomorphic to
the  direct sum of two
copies of matrix algebra $M_{\ell}(\C)$. Then we shows the quotient algebra
$A(\calh^+)/\cali$ is commutative and is generated by the
elements
$\omega_a, J_a$ and $\Lambda_{ab}$ (see Section 5). It is fair to say that
we do not determine the algebra structure of $A(\calh^+)$ 
completely in terms of generators
and relations. But the relations among generators of
$A(\calh^+)$ found in this paper are
good enough to classify all irreducible modules for $A(\calh^+)$ and
for the VOA $\calh^+$.

We organize the paper as follows. In Section 2, we review definitions and
states properties of the VOA $\calh^+$. The list of inequivalent
irreducible modules is given here. We explain the notion of Zhu's algebras
in Section 3 and prove some formulas which we need later. Section 4 is
devoted to find a finite set of generators for $A(\calh^+).$ In Section 5, we
introduce the elements $E_{ab}^u, E_{ab}^t$ and $\Lambda_{ab}$ as well as
$\omega_a, J_a$ which form a ``nice" generating set of Zhu's algebra.
It will be shown that the elements $E_{ab}^u$ and $E_{ab}^t$ forms the
matrix algebra $M_\ell(\C)$ respectively. We derive more relations
among the generators in Section 6. The evaluation method developed in
\cite{DN1} and \cite{DN2} is fully used for this aim. We show these
relations are enough to classify all the irreducible modules for
$A(\calh^+)$ and then for the VOA $\calh^+$.

The core of this work was done while the second author was visiting University
of California at Santa Cruz in January, 1999. K.~N thanks 
Professor Mason for the hospitality during the stay.

Throughout the paper $\N$ is the set of nonnegative integers and $\Z_+$
is the set of positive integers.

\section{Preliminaries}
\label{section:2}
This section is divided into two subsections. In the first subsection we
recall various notions of (twisted) modules for a vertex operator algebra
$V$. In the second subsection we discuss the construction of 
the vertex operator algebra $\calh$ and its (twisted) modules. 

\subsection{Modules}\label{subsection:2.1}
Let $V=\oplus_{n\in\Z}V_n$ be a vertex operator algebra
(cf.~\cite{B}, \cite{FLM}) and $g$ be an automorphism of $V$ of finite
order $T$. Then $g$ preserves each weight space $V_n$ and we decompose
$V$ into eigenspaces with respect to the action of $g$ as
$V=\bigoplus_{r\in \Z/T\Z}V^r$ where $V^r=\{v\in V|gv=e^{-2\pi ir/T}v\}$.

An \textit{admissible $g$-twisted $V$-module} (cf.~[DLM2])  
\[
M=\sum_{n=0}^{\infty}M(\frac{n}{T})
\]
is a $\frac{1}{T}\N$-graded vector space with the \textit{top level} 
$M(0)\neq 0$ equipped with a linear map
\[
\begin{array}{ccl}
V&\longrightarrow &(\End\,M)\{z\}\\
v&\longmapsto& Y_M(v,z)=\sum_{n\in\Q}v_nz^{-n-1},\ \ \ (v_n\in
\End\,M)
\end{array}
\]
which satisfies the following conditions for
$0\leq r\leq T-1,$ $u\in V^r$, $v\in V,$ 
$w\in M$:

\noindent
(A1) $Y_M(u,z)=\sum_{n\in r/T+\Z}u_nz^{-n-1}$,
i.e., $u_n = 0$ if $n\notin r/T+\Z$.

\noindent
(A2) There exists an integer $N$ such that $u_nw=0$ for all
$n>r/T+N$.

\noindent
(A3) $Y_M(\unit,z)=\id_M$ where $\id_M$ is the identity map on $M$.

\noindent
(A4) Jacobi identity
\[\label{eqn:2.1}
\begin{split}
z^{-1}_0\delta\left(\frac{z_1-z_2}{z_0}\right)
&Y_M(u,z_1)Y_M(v,z_2)-z^{-1}_0\delta\left(\frac{z_2-z_1}{-z_0}\right)
Y_M(v,z_2)Y_M(u,z_1)\\
&=z_2^{-1}\left(\frac{z_1-z_0}{z_2}\right)^{-r/T}
\delta\left(\frac{z_1-z_0}{z_2}\right)
Y_M(Y(u,z_0)v,z_2),
\end{split}
\]
where $\delta(z)=\sum_{n\in\Z}z^n$ and all binomial expressions are to be 
expanded in nonnegative
integral powers of the second variable. (One can find elementary
properties of the $\delta$-function in \cite{FLM}.)

\noindent
(A5) If $u$ is homogeneous,
\[
u_mM(n)\subset M(n+\wt{u}-m-1).
\]

\vskip 1ex

If $g=\id_V$, this reduces to the definition of an admissible $V$-module
(cf. [DLM1]).

A $g$-{\em twisted $V$-module} is
an admissible $g$-twisted $V$-module $M$ such that $L(0)$ is semisimple;
\[
M=\coprod_{\lambda \in{\C}}M_{\lambda},\quad  M_{\l}=\{w\in M|L(0)w=\l w\}
\]
and $\dim M_{\l}$ is finite, and for fixed $\l,$ $M_{{n/T}+\l}=0$
for all small enough integers $n$. Again if $g=\id_V$  we get the
definition of an ordinary
$V$-module.

\subsection{Vertex operator algebras $\calh$ and $\calh^+$}
\label{subsection:2.2}

Following \cite{FLM} we discuss the construction of vertex operator algebra
$\calh$ and its (twisted) modules. The vertex operator subalgebra
$\calh^+$ is defined and the list of known irreducible modules for
$\calh^+$ is presented. 

Let ${\frakh}$ be an $\ell$-dimensional vector space
with a nondegenerate symmetric bilinear form 
$\<\,,\,\>$ and $\hat {\frakh}={\frakh}\otimes 
\C[t,t^{-1}]\oplus \C\,K$ be the corresponding affine Lie algebra viewing
$\frakh$ as an abelian Lie algebra. Let $\l\in \frakh$ and consider the
induced $\hat {\frakh}$-module
\[
%\label{eqn:2.8}
M(1,\l)= 
U(\hat{\mathfrak h})\otimes_{U({\mathfrak h}\otimes{{\C}}[t]
\oplus{{\C}}K)}{{\C}}\simeq S(\mathfrak h\otimes t^{-1}\C[t^{-1}])\ \ \
(\mbox{linearly})
\]
where ${\mathfrak h}\otimes t{{\C}}[t]$ acts trivially on $\C,$
${\mathfrak h}$ acts as $\<\alpha,\l\>$ for $\alpha\in{\mathfrak h}$ and
$K$ acts as $1$. For $\alpha\in{\mathfrak h}$ and $n\in {\Z}$,  we write
$\alpha(n)$ for the operator $\alpha\otimes t^n$ acting on
$M(1,\lambda)$ and set
\[
\alpha(z)=\sum_{n\in{{\Z}}}\alpha(n)z^{-n-1}.
\]

Among $M(1,\lambda),(\lambda \in \mathfrak{h})$, the space 
$\calh = M(1,0)$ is specially interesting as it has a vertex
operator algebra structure as explained below. ($M(1,0)$ is denoted 
by $M(1)$ in [FLM].) 
We set $\unit = 1\otimes 1$.
For $\alpha_1,\dots,\alpha_k\in {\mathfrak h},\  (n_1, \dots, n_k
\in{\Z_+})$ and $v=\alpha_1(-n_1)\cdots\alpha_k(-n_k)\unit\in \calh$,
we define a vertex operator, acting on $M(1,\lambda),$ 
 corresponding to $v$ by
\begin{equation}\label{eqn:2.9}
Y(v,z)=\NO[
\partial^{(n_1-1)}\alpha_1(z)][\partial^{(n_2-1)}\alpha_2(z)]\cdots
[\partial^{(n_k-1)}\alpha_k(z)]\NO,\quad
\partial^{(n)} = \frac{1}{n!}\left(\frac{d}{dz}\right)^n
\end{equation}
where a normal ordering procedure indicated by open colons signifies
that the expression between two open colons is to be reordered if
necessary so that all the operators
$\alpha(n)$ $(\alpha\in{\mathfrak h},\  n<0)$  are to be placed to
the left of all the operators $\alpha(n),\,(n\ge 0)$ before the expression 
is evaluated. We extend $Y$ to all $v\in V$ by linearity. 
Let $\{h_1,...h_{\ell}\}$ be an orthonormal basis of $\frakh$ and set 
$\omega=\frac{1}{2}\sum_{i=1}^{\ell}h_i(-1)^2\unit$. The following
theorem is well known (cf.~\cite{FLM}).

\begin{theorem}\label{theorem:2.4}
The space $\calh=(\calh,Y,{\bf 1},\omega)$ is a simple vertex 
operator algebra with a vacuum $\unit$ and a Virasoro element $\omega$ and
$M(1,\l)=(M(1,\l),Y)$ for $\l\in
\frakh$  are inequivalent $\calh$-modules. Moreover, 
any irreducible $\calh$-module is isomorphic to a module 
$M(1,\l)=(M(1,\l),Y)$ for for some $\l\in
\frakh.$  
\end{theorem}

Now, we define an automorphism $\theta$ of $\calh$ by 
\[
\theta(\a_1(-n_1)\cdots
\a_k(-n_k)\unit)=(-1)^k\a_1(-n_1)\cdots\a_k(-n_k)\unit.
\]
Then the $\theta$-fixed point set $\calh^+$ of $\calh$ is a simple vertex
operator subalgebra and the $-1$-eigenspace $\calh^-$ is an irreducible
$\calh^+$-module: See Theorem~2 of \cite{DM2}.

Following \cite{DM1} we define another $\calh$-module from 
$M(1,\lambda)$; $\theta\circ
M(1,\l)=(\theta\circ M(1,\l), Y^{\theta})$ where $Y^{\theta}(v,z) =Y(\theta
(v),z)$. Then
$\theta\circ M(1,\l)$ is also an irreducible $\calh$-module isomorphic to
$M(1,-\l)$. The following proposition is a direct consequence of Theorem
6.1 of
\cite{DM2}.
\begin{proposition}\label{proposition:2.5} 
If $\l\ne 0$ then $M(1,\l)$ is an irreducible $\calh^+$-module, and
$M(1,\l)$ and $M(1,-\lambda)$ are isomorphic.
\end{proposition}

Next we turn our attention to the $\theta$-twisted $\calh$-modules
(cf.~[FLM]). The twisted affine algebra for $\frakh$ is defined to be
$\hat{\frakh}[-1]=
\sum_{n\in\Z}{\mathfrak h}\otimes t^{1/2+n}\oplus \C\,K$. Its canonical 
irreducible module is 
$$\calh(\theta)=U(\hat {\mathfrak h}[-1])\otimes_{U({\mathfrak h}\otimes
t^{1/2}\C[t]\oplus \C 
K)}\C\simeq S({\mathfrak h}\otimes t^{-1/2}\C[t^{-1/2}])$$
where ${\mathfrak h}\otimes t^{1/2}\C[t]$ acts trivially on $\C$ and
$K$ acts as $1$. As before we can define an involution on
$\calh(\theta)$ also denoted
by $\theta$  
\[
\theta
(\a_1(-m_1)\cdots \a_k(-m_k)\unit) =(-1)^k\a_1(-m_1)\cdots \a_k(-m_k)\unit
\]
where $\a_i\in {\mathfrak h}, m_i\in 1/2+\N$ and
$\alpha(n)=\alpha\otimes t^n$. We denote the $\pm 1$-eigenspace of
$\calh(\theta)$ under $\theta$  by $\calh(\theta)^{\pm}.$

Let $v=\a_1(-n_1)\cdots \a_k(-n_k)\unit\in
\calh,\,(n_1,n_2,\dots,n_k\in\Z_+)$. We first introduce an operator 
\[
%\label{eqn:2.6*}
W_{\theta}(v,z)=\NO
[\partial^{(n_1-1)}\alpha_1(z)][\partial^{(n_2-1)}\alpha_2(z)]
\cdots[\partial^{(n_k-1)}
\alpha_k(z)]\NO,
\]
where the right hand side is an operator on $\calh(\theta)$, namely,
\[
\alpha(z)=\sum_{n\in\frac{1}{2}+\Z}\alpha(n)z^{-n-1}
\]
and where normal ordering procedure is as before. We extend
this to all $v\in \calh$ by linearity. Define constants $c_{mn}\in \Q$
for $m, n\ge0$ by the formula
\[
\sum_{m,n\ge 0}c_{mn}x^my^n=-{\rm log}\left(\frac
{(1+x)^{1/2}+(1+y)^{1/2}}{2}\right)
\]    
and set
\[
\Delta_z=\sum_{m,n\ge 0}\displaystyle{
\sum^\ell_{i=1}}
c_{mn}h_i(m)h_i(n)z^{-m-n}.
\]
The \textit{twisted vertex operator}
$Y_{\theta}(v,z)$ for $v\in \calh$ is defined by
\[
Y_{\theta}(v,z)=W_{\theta}(e^{\Delta_z}v,z).
\]
 
Then we have:
\begin{theorem}\label{theorem:2.6} {\rm (i)} $(\calh(\theta),Y_{\theta})$ is 
an irreducible $\theta$-twisted $\calh$-module.

\noindent
{\rm (ii)} $\calh(\theta)^{\pm}$ are irreducible $\calh^+$-modules.
\end{theorem}
\begin{proof}
Part (i) is a result of Chapter 9 of \cite{FLM} and part (ii) follows from
Theorem 5.5 of \cite{DL}.
\end{proof}

\section{Zhu's algebra}
We review the definition of Zhu's algebra $A(V)$ associated to a vertex
operator algebra $V$ and related results from [Z] and [DLM2].
We also give several frequently used formulas in $A(\calh^+).$ 

\subsection{The definition of Zhu's algebra}
Let us recall a vertex operator algebra is $\Z$-graded;
\[
V = \oplus_{n\in\Z}V_n.
\]
Each $v\in V_n$ is called a homogeneous element of $V$ with weight $n$,
which we denote $n = \wt{v}$. Whenever we write $\wt{v}$, the element $v$
is assumed to be homogeneous. In order to introduce Zhu's algebra, we need
to define two binary operations $*$ and $\circ$ on $V$. 

For $u\in V$ homogeneous and $v\in V$, we define
\begin{align}
u*v &= \Res{z}\left(\frac{(1+z)^{\wt{u}}}{z}Y(u,z)v\right) =
\sum_{i=0}^\infty\binom{\wt{u}}{i}u_{i-1}v,\label{eqn:3.1}\\
u\circ v &= \Res{z}\left(\frac{(1+z)^{\wt{u}}}{z^2}Y(u,z)v\right) =
\sum_{i=0}^\infty\binom{\wt{u}}{i}u_{i-2}v\label{eqn:3.2}
\end{align}
and extend both (\ref{eqn:3.1}) and (\ref{eqn:3.2}) to linear operations on
$V.$ Define $O(V)$ to be the linear space spanned by all
$u\circ v$ for $u,v\in V$. Then  the $A(V)$ is defined to be the
quotient space $V/O(V)$.

For $u\in V$, we define $o(u)$ the weight zero component operator of $u$ on
any admissible modules. Then $o(u) = u_{\wt{u}-1}$ if $u$ is homogeneous.

The following theorem is essentially due to Zhu \cite{Z} (also see
\cite{DLM2}).

\begin{theorem}\label{theorem:3.1.1}
{\rm (i)}
The bilinear map $*$ induces an associative algebra structure
on $A(V)$ with the identity $\unit + O(V)$. Moreover $\omega+ O(V)$ is a
central element of $A(V)$.

\noindent
{\rm (ii)} 
The map $u\longmapsto o(u)$ gives a representation of $A(V)$ on
$M(0)$ for any admissible $V$-module $M$. Moreover, if any admissible
module is completely reducible, then  $A(V)$ is a finite dimensional
semisimple algebra.

\noindent
{\rm (iii)} 
The map $M\longmapsto M(0)$ gives a bijection between the set of
equivalence classes of irreducible admissible $V$-modules and the set of
equivalence classes of simple $A(V)$-modules.
\end{theorem}

For convenience, we write $[u] = u + O(V)\in A(V)$. We denote $u\sim v$
for $u,v\in V$ if $[u] = [v]$. This induces an equivalence relation on
$\End\, V$ such that for $f,g \in \End\,V$, $f\sim g$ if and only if
$fu\sim gu$ for all $u\in V$.

The following proposition is useful later (see \cite{Z} for details).

\begin{proposition}\label{proposition:3.1.2}
{\rm (i)} Let $u\in V$ be homogeneous and $v\in V$, then for $n\in \N$
\[
\Res{z}\left(\frac{(1+z)^{\wt{u}}}{z^{n+2}}Y(u,z)v\right) =
\sum_{i=0}^\infty\binom{\wt{u}}{i}u_{i-n-2}v\in O(V).
\]

\noindent
{\rm(ii)}
Let $v\in V$ be homogeneous and $u\in V$, then
\[
u*v\sim \Res{z}\left(\frac{(1+z)^{\wt{v}-1}}{z}Y(v,z)u\right) =
\sum_{i=0}^\infty\binom{\wt{v}-1}{i}v_{i-1}u.
\]

\noindent
{\rm (iii)}
Let $u,v\in V$ be homogeneous, then
\[
u*v -v*u \sim \Res{z}\left((1+z)^{\wt{u}-1}Y(u,z)v\right) =
\sum_{i=0}^\infty \binom{\wt{u}-1}{i}u_{i}v.
\]

\noindent
{\rm (iv)} For any $u\in V,$ $L(-1)u+L(0)u\in O(V)$ where $Y(\omega,z)
= \sum_{n\in\Z}L(n)z^{-n-2}$. 
\end{proposition}

\subsection{Some formulas in $A(\calh^+)$}
We prove some formulas in $A(\calh^+)$.

Recall that $\{h_a\in \frakh\,|\,a = 1,2,\dots, \ell\,\}$ is an orthonormal
basis. We set
\[
\omega_a = \frac{1}{2}h_a(-1)^2\unit \quad \text{and} \quad Y(\omega_a,z)
= \sum_{n\in \Z}L_a(n) z^{-n-2}
\]
for $a = 1,2,\dots,\ell$.

\begin{proposition}\label{proposition:3.2.1}
For all $a,\,(1\leq a\leq \ell)$,

\noindent
{\rm (i)}
$ L_a(-n-3) + 2L_a(-n-2) + L_a(-n-1)\sim 0,\, (n\geq 0)$.

\noindent
{\rm (ii)}
$u*\omega_a\sim (L_a(-2) + L_a(-1))u$ for all $u\in \calh^+$.

\noindent
{\rm (iii)}
$\omega_a*u-u*\omega_a\sim (L_a(-1)+L_a(0))u$ for all $u\in\calh^+$.
\end{proposition}
\begin{proof}
Using Proposition \ref{proposition:3.1.2} (i) for $\omega_a$ and for any
$u\in V$ proves
\[
\Res{z}\left( \frac{(1+z)^2}{z^{2+n}}Y(\omega_a,z)u\right) =
\bigl(L_a(-n-3) + 2L_a(-n-2) + L_a(-n-1)\bigr)u\in O(V).
\]  
Proposition \ref{proposition:3.1.2} (ii) with $v = \omega_a$ shows
\[
u*\omega_a\sim \Res{z}\left(\frac{1+z}{z}Y(\omega_a,z)u\right)
=(L_a(-2)+L_a(-1))u.
\]
Finally by Proposition \ref{proposition:3.1.2} (iii), we have
\[
\omega_a*u-u*\omega_a \sim \Res{z}\left((1+z)Y(\omega_a,z)u\right) =
(L_a(-1)+L_a(0))u.
\]
\end{proof}

\begin{remark}\label{note:3.2.2}
The Virasoro element $\omega$ in $\calh^+$ is the sum of
$\omega_a$'s, i.e., $\omega = \sum_{a=1}^\ell \omega_a$. We know
$L(-1)+L(0)\sim 0$. But it is less obvious that $L_a(-1)+L_a(0)\nsim 0$. 
In order to see this, we fix  $\{h_1(-1)\unit,\dots,h_{\ell}(-1)\unit\}$
as a basis of the top level of $\calh^{-}$ and set 
$S_{ab} = h_a(-1)h_b(-1)\unit.$ Then $S_{ab}$ acts on the top level of
the module $\calh^-$ as $E_{ab}+E_{ba}$ with respect to the basis
where $E_{ab}$ is the matrix element sending $h_c(-1)$ to 
$\delta_{c,a}h_b(-1),$ and $\omega_a$ acts as
$E_{aa}$: See Table 1 in Subsection \ref{subsection:4.3}. Therefore,
$\omega_a*S_{ab} - S_{ab}*\omega_a$ acts as
$E_{ab}$, which means $(L_a(-1)+L_a(0))S_{ab}\nsim 0$ by Proposition
\ref{proposition:3.2.1} (iii). In particular, the algebra $A(\calh^+)$ is
not commutative.
\end{remark}

\section{A finite set of generators for  $A(\calh^+)$}
\label{section:4}

In this section we prove that the algebra $A(\calh^+)$ is finitely generated.
The main
idea is to consider a vertex operator subalgebra $\calw$ of $\calh^+$
such that  $\calh^+$ is a finite direct sum of irreducible
modules for $\calw.$ We use a result
from [DN1] to show that the image $A[\calw]$ of $\calw$ in
$A(\calh^+)$ is finitely generated. Finally we determine a finite set
of generators for the image of each irreducible $\calw$-submodule
in $A(\calh^+)$ as a left or right $A[\calw]$-module. 

For convenience, we sometimes identify an element $u$ in a vertex
operator algebra $V$ with its image $[u]=u+O(V)$ in $A(V)$ if there
is no confusion arising. 

\subsection{A generating set of $A(\calh^+)$}
\label{subsection:4.1}

Let $\mathcal{H}_a $ be the vertex operator subalgebra algebra (with
Virasoro element $\omega_a$) associated to
the 1-dimensional vector space $\C h_a$. 
Then the automorphism $\theta$ of
$\calh$ induces an automorphism of $\calh_a$ denoted by the same symbol
$\theta$ and $\mathcal{H}_a$ decomposes into direct sum of 
the $\pm 1$-eigenspaces of $\theta$;
\[
\mathcal{H}_a = \mathcal{H}_a^+ \oplus \mathcal{H}_a ^-.
\]
As a result, $\mathcal{H}$ decomposes into
\[
%\label{eqn:3}
\mathcal{H} = \bigoplus_{\alpha\subset \{1,...,\ell\}}{\calw}^{\alpha} 
\]
where 
\[
{\calw}^{\alpha}= 
\mathcal{H}_1^{\varepsilon_1}\otimes \cdots \otimes
\mathcal{H}_\ell^{\varepsilon_\ell}
\]
such that $\varepsilon_i=-$ if $i\in\alpha$ and  $\varepsilon_i=+$ if 
$i\not\in\alpha.$ 

For convenience, we also write ${\calw}={\calw}^{\emptyset}=
\calh_1^+\otimes\cdots\otimes\calh_\ell^+$. Let $P$ be the collection
of all subsets of $\{1,...,\ell\}$ with even cardinalities. 
Then 
\[
%\label{eqn:3'}
\mathcal{H}^+ = \oplus_{\alpha\in P}{\calw}^{\alpha}, 
\]
the ${\calw}$ is a vertex operator subalgebra of $\calh^+$ and each 
${\calw}^{\alpha}$ is an irreducible ${\calw}$-submodule of $\calh^+.$
For a subspace $W$ of $\calh^+$ we denote the image of $W$ in
$A(\calh^+)$ by $A[W].$ That is, $A[W]=\{w+O(\calh^+)|w\in W\}.$ 
Then
\[
%\label{eqn:3''}
A(\calh^+)=\oplus_{\alpha\in P}A[\calw^{\alpha}],
\]
$A[\calw]$ is a subalgebra of $A(\calh^+)$ and each $A[\calw^{\alpha}]$
is an $A[\calw]$-bimodule of $A(\calh^+).$ 

The main purpose in this subsection is to find a set of generators
of each $A[\calw^{\alpha}]$ as a left or right $A[\calw]$-module. 

We first find a set of generators for the algebra $A[\calw].$
We need the following Lemma from \cite{DMZ}.

\begin{lemma}\label{lemma:4.1.1}
Let $V_1,\dots, V_n$ be vertex operator algebras. Then the linear map
\[
F:[v_1]\otimes \cdots\otimes [v_n]\longmapsto [v_1\otimes \cdots \otimes
v_n]
\]
from $A(V_1)\otimes \cdots \otimes A(V_n)$ to $A(V_1\otimes \cdots
\otimes V_n)$ is an isomorphism of associative algebras.
\end{lemma}

Let $J_a = h_a(-1)^4\unit -2h_a(-3)h_a(-1)\unit +
\frac{3}{2}h_a(-2)^2\unit$. 
Then the following lemma was shown in \cite{DN1}.

\begin{lemma}\label{lemma:4.1.2a}
For $a = 1,2,\dots,\ell$, $A(\calh_a^+)$ is generated by $\omega_a$ and
$J_a$. 
\end{lemma}

\begin{corollary}\label{corollary:4.1.2}
The $A[\mathcal{W}]$ is generated by $\omega_a$ and $J_a$ for $a =
1,\dots,\ell$. In particular, $A[\calw]$ is a commutative 
subalgebra of $A(\calh^+).$ 
\end{corollary}
\begin{proof}
By Lemma \ref{lemma:4.1.1}  the composition of the linear map
\[
\begin{array}{cccc}
F:&A(\calh_1^+)\otimes\cdots\otimes
A(\calh_\ell^+)&\longrightarrow&A(\calh_1^+\otimes\cdots\otimes
\calh_\ell^+) = \mathcal{W}/O(\mathcal{W})\\
&[v_1]\otimes \cdots\otimes [v_n]&\longmapsto
&[v_1\otimes \cdots \otimes v_n]
\end{array}
\]
and the canonical map
$\iota:A(\calh_1^+\otimes\cdots\otimes \calh_\ell^+)\longrightarrow
A[\calw]
$ which sends $u+O(\calw)$ to $u+O(\calh^+)$ 
is a surjective homomorphism. Then Lemma \ref{lemma:4.1.2a} shows that
$A[\mathcal{W}]$ is generated by $\omega_a$ and $J_a$ for $a=
1,2,\dots,\ell$.
\end{proof}

Next, we consider $A[\calw^{\alpha}]$ with $\alpha=\{a,b\}$
for distinct $a,b,\,(1\leq a,b\leq \ell)$. 
For positive integers $m,n\in\Z_+$, we set
$S_{ab}(m,n) = h_a(-m)h_b(-n)\unit$.

%\begin{lemma}\label{lemma:4.1.3}
%For $a,b, \,(1\leq a,b\leq \ell)$ distinct and integers $m,n\geq 1$,
%the element $S_{ab}(m,n)$ is a linear combination of elements
%\[
%[\,\omega_a,[\,\omega_a,\cdots, [\,\omega_a,S_{ab}(1,1)\,],\cdots,]]
%\]
%in $A[\calw^{\alpha}].$
%\end{lemma}
%\begin{proof}
%We prove the lemma by induction on the weight $N=m+n$ of $S_{ab}(m,n)$. 
%Clearly the
%minimal weight is $2$ and $S_{ab}(m,n) = h_a(-1)h_b(-1)\unit = S_{ab}(1,1)$
%in this case.

%Now let $N>2$ and suppose that the lemma
%is true for all  $S_{ab}(m,n)$ subject to $m+n<N$. For $S_{ab}(m,n)$ with
% $m+n = N$, one can assume that $m\geq 2$ without loss of generality.
%Then, since 
%\begin{align*}
%&L_a(-1)S_{ab}(m-1,n) = (m-1)S_{ab}(m,n)\quad\text{and}\\
%&L_a(0)S_{ab}(m-1,n) = (m-1)S_{ab}(m-1,n),
%\end{align*}
%Proposition \ref{proposition:3.2.1} shows
%\begin{align*}
%\omega_a*S_{ab}(m-1,n) &- S_{ab}(m-1,n) *\omega_a\\
%&\sim (L_a(-1)+L_a(0))S_{ab}(m-1,n) \\
%&= (m-1)S_{ab}(m,n) + (m-1)S_{ab}(m-1,n).
%\end{align*}
%By the induction hypothesis, we see that $S_{ab}(m,n)$ is a linear
%combination of the prescribed elements.
%\end{proof}

%\begin{remark}
%As one can see in the proof of Lemma \ref{lemma:4.1.3}, the element
%$S_{ab}(m,n)$ is also a linear combination of the elements
%\[
%[\,\omega_b,[\,\omega_b,\cdots, [\,\omega_b,S_{ab}(1,1)\,],\cdots,]].
%\]

%Alternately one can show this by noting that the Virasoro element $\omega$ 
%is a
%central element of $A(\calh^+)$ and in particular, $[\omega, S_{ab}(1,1)\,]
%= 0$, which is equivalent to $[\omega_a,S_{ab}(1,1)\,] =
%-[\omega_b,S_{ab}(1,1)\,]$.
%\end{remark}

Let $\cals_{ab}$ be the
linear subspace spanned by $S_{ab}(m,n)+O(\calh^+),\,(m,n\in\Z_+).$ 

\begin{lemma}\label{lemma:4.1.5} We have $A[\calw^{\alpha}]=
\cals_{ab}A[\calw]$ for $\alpha=\{a,b\}.$ That is, 
as a right $A[\calw]$-module, $A[\calw^{\alpha}]$ is generated by
$\cals_{ab}.$
\end{lemma}
\begin{proof} 
Note that $\calw^{\alpha}$ is spanned by 
$S_{ab}(m,n;u)=h_a(-m)h_b(-n)u$ where $u\in\mathcal{W}$ and $m,n\geq 1$
are integers. We can assume that $u$ is 
a monomial $u = h_{a_1}(-n_{1})\cdots  h_{a_r}(-n_{r})\unit$ in
$\mathcal{W}.$ For such $u$  we define the length $\ell(u)$ to be $r.$  
As usual
we define $\ell(\unit) = 0$.

We prove  by induction on the length of
monomial $u\in\mathcal{W}$ that $S_{ab}(m,n;u)+O(\calh^+)\in
\cals_{ab}A[\calw].$ 
If $\ell(u) = 0$, then 
$S_{ab}(m,n;u) = S_{ab}(m,n)$ and
it is clear.

Suppose the lemma is true for all monomials with lengths strictly less
than $N$. Let $u\in \mathcal{W}$ with the length $\ell(u)=N$. Recall that 
the component operator $S_{ab}(m,n)_i$ of $Y(S_{ab}(m,n),z)$ is defined
by
\[
Y(S_{ab}(m,n),z) = \sum_{i\in\Z}S_{ab}(m,n)_iz^{-i-1}.
\]
It follows from (\ref{eqn:2.9}) that 
\[
S_{ab}(m,n)_{i-1} = \sum_{\substack{
r+s= i-m-n\\
r\geq 0\,\text{or}\,r\leq -m\\
s\geq 0\,\text{or}\,s\leq -n}}d_{rm}d_{sn}\NO h_a(r)h_b(s)\NO,\quad d_{rm}
= \binom{-r-1}{m-1}
\]
for $i\in\Z.$
Now we compute the $*$ product
\[
S_{ab}(m,n)*u = \sum_{i=0}^{m+n} \binom{m+n}{i} S_{ab}(m,n)_{i-1}u.
\]
Note that if $r\leq -m$ and $s\leq -n$, then $r+s\leq -m-n$. Therefore,
if $i=0$, then
\[
S_{ab}(m,n)_{-1}u = S_{ab}(m,n;u) + \sum_{\substack{
r+s= -m-n\\
r\geq 0\,\text{or}\,s\geq0}}
d_{rm}d_{sn}\NO h_a(r)h_b(s)\NO u,
\]
and if $i>0$, then we see either $r\geq 0$ or $s\geq 0$, namely,
\[
S_{ab}(m,n)_{i-1}u =  \sum_{\substack{
r+s= i-m-n\\
r\geq 0\,\text{or}\,s\geq0}}
d_{rm}d_{sn}\NO h_a(r)h_b(s)\NO u.
\]
Hence we have $S_{ab}(m,n)*u = S_{ab}(m,n;u) + w$ where $
w =
\sum_{r,s\geq 1} h_a(-r)h_b(-s)u^{r,s}$ and where $u^{r,s}\in \mathcal{W},
\,\ell(u^{r,s})<N$. Then by the induction hypothesis $w+O(\calh^+)\in 
\cals_{ab}A[\calw].$
As a result, we have $S_{ab}(m,n;u)+O(\calh^+)\in \cals_{ab}A[\calw]$.
\end{proof}

The following proposition is similar to Lemma \ref{lemma:4.1.5}.
\begin{proposition}\label{proposition:4.1.6} We have
$A[\calw^{\alpha}]=A[\calw] \cals_{ab}$ for $\alpha=\{a,b\}.$ That is, 
as a left $A[\calw]$-module, $A[\calw^{\alpha}]$ is generated by
$\cals_{ab}.$
\end{proposition}
\begin{proof}
It is enough to prove $ \cals_{ab}A[\calw]\subset A[\calw]  \cals_{ab}$.
Recall the
vector $S_{ab}(m,n;u)$ from the proof of Lemma \ref{lemma:4.1.5}.
We also use
induction on the length of the monomial
$u$. If
$\ell(u) = 0$, it is clear. Let $N$ be a positive integer and suppose that
the claim is true for all monomials $u\in \mathcal{W}$ with lengths
strictly less
than $N$. Now, consider
$S_{ab}(m,n;u)$ for $u\in \mathcal{W}$ with $\ell(u) = N$. Proposition
\ref{proposition:3.1.2} (iii) shows 
\[
\begin{split}
&S_{ab}(m,n)*u-u*S_{ab}(m,n)\\
&\sim
\Res{z}\left((1+z)^{m+n-1}Y(S_{ab}(m,n),z)u\right) =
\sum_{i=0}^{m+n-1}\binom{m+n-1}{i}S_{ab}(m,n)_{i}u.
\end{split}
\]
{}From the proof of Lemma \ref{lemma:4.1.5} we see that
\[
S_{ab}(m,n)_{i}u = \sum_{r,s\geq
1}h_a(-r)h_b(-s)u^{r,s}\quad\text{for}\quad i\geq 0
\]
where $u^{r,s}\in \mathcal{W}$ and $\ell(u^{r,s})<N$. Thus 
each $S_{ab}(m,n)_{i}u+O(\calh^+)$ is a linear combination
of $S_{ab}(s,t)*v+O(\calh^+)$'s where $s,t>0$ and where $v\in \calw$ 
are monomials with lengths less than $N.$ 
By the induction
hypothesis each $S_{ab}(m,n)_{i}u+O(\calh^+)$ lies in $A[\calw]\cals_{ab}.$
Thus $S_{ab}(m,n)*u -u*S_{ab}(m,n)+O(\calh^+)\in A[\calw]\cals_{ab}$ and  
$S_{ab}(m,n)*u+O(\calh^+)\in A[\calw]\cals_{ab}.$ 
\end{proof}

We now turn our attention to $A[\calw^{\alpha}]$ for general $\alpha.$
For this purpose we consider the elements of type
\[
S_{abcd}(m,n,r,s) = h_a(-m)h_b(-n)h_c(-r)h_d(-s)\unit
\]
where $m,n,r,s\in\Z_{+}$ and $a,b,c,d$ are distinct.

\begin{lemma}\label{lemma:4.1.6}
For any $m,n,r,s\in\Z_{+}$,
\[
S_{abcd}(m,n,r,s)+O(\calh^+)=(-1)^{m+n+r+s}S_{abcd}(1,1,1,1)+O(\calh^+)
\]
\end{lemma}
\begin{proof} Recall the definition of the circle operation
\[
S_{ab}(m,n)\circ S_{cd}(r,s) = \sum_{k=0}^{m+n}
\binom{m+n}{k} S_{ab}(m,n)_{k-2}S_{cd}(r,s).
\]
Also recall that
\[
S_{ab}(m,n)_{k-2} = \sum_{\substack{
i+j= k-1-m-n\\
i\geq 0\,\text{or}\,i\leq -m\\
j\geq 0\,\text{or}\,j\leq -n}}
d_{im}d_{jn}\NO h_a(i)h_b(j)\NO,\quad d_{im} = \binom{-i-1}{m-1}
\]
and note that if $k\geq 2$ then either $i\geq 0$ or $j\geq 0$ in the sum.
This immediately gives $S_{ab}(m,n)_{k-2}S_{cd}(r,s) = 0$ for $k\geq 2$. 
Thus we
have $S_{ab}(m,n)\circ S_{cd}(r,s) = S_{ab}(m,n)_{-2}S_{cd}(r,s)+
(m+n)S_{ab}(m,n)_{-1}S_{cd}(r,s).$ It is easy to see that
\begin{align*}
S_{ab}(m,n)_{-2}S_{cd}(r,s)& = mS_{abcd}(m+1,n,r,s) + n
S_{abcd}(m,n+1,r,s),\quad\text{and}\\
S_{ab}(m,n)_{-1}S_{cd}(r,s)&= S_{abcd}(m,n,r,s).
\end{align*}
So
\begin{equation}\label{eqn:4'}
mS_{abcd}(m+1,n,r,s) + nS_{abcd}(m,n+1,r,s) +(m+n)S_{abcd}(m,n,r,s)\sim 0.
\end{equation}

Similarly, when considering
$S_{ac}(m,r)\circ S_{bd}(n,s)$ and $S_{bc}(n,r)\circ S_{ad}(m,s)$
respectively, we obtain
\begin{equation}\label{eqn:4}
mS_{abcd}(m+1,n,r,s) + rS_{abcd}(m,n,r+1,s) +
(m+r)S_{abcd}(m,n,r,s)\sim 0
\end{equation}
and 
\begin{equation}\label{eqn:4''}
nS_{abcd}(m,n+1,r,s) + rS_{abcd}(m,n,r+1,s) +
(n+r)S_{abcd}(m,n,r,s)\sim 0.
\end{equation}
Add (\ref{eqn:4'}) and (\ref{eqn:4}) together and use (\ref{eqn:4''})
to yield
$$S_{abcd}(m+1,n,r,s) +
S_{abcd}(m,n,r,s)\sim 0.$$
Consequently we have $S_{abcd}(m,n,r,s)\sim
(-1)^{m-1}S_{abcd}(1,n,r,s)$.
Since $S_{abcd}(m,n,r,s)$ is invariant under the permutations of
$\{(a,m), (b,n), (c,r), (d,s)\}$, we
can apply the same result to indices $b,c$ and $d$ and finish the 
proof of the lemma.
\end{proof}

We denote $S_{ab} = S_{ab}(1,1)$ for short.

\begin{remark}\label{remark:4.1.9}
We see from Lemma \ref{lemma:4.1.6} that
\[
S_{abcd}(m,n,r,s) \sim (-1)^{m+n+r+s}S_{ab}*S_{cd}.
\]
as $S_{abcd}(1,1,1,1) = S_{ab}*S_{cd}$.
\end{remark}

Now let $\alpha=\{a_1,...,a_{2k}\}$  be an even subset of
$\{\,1,2,\dots,\ell\,\}.$ Let $\alpha = \alpha_1\cup
\alpha_2\cup\cdots
\cup
\alpha_k$ be a disjoint union of subsets $\alpha_i$ such that 
$|\alpha_i| = 2.$ 
Set $S_{\alpha}=S_{\alpha_1}*S_{\alpha_2}*\cdots *S_{\alpha_k}$ where
where $S_\alpha = S_{ab}$ for $\alpha = \{a,b\}$. 
  Clearly
$S_\alpha$ is independent of a choice of decomposition $\alpha =
\alpha_1\cup \alpha_2\cup\cdots \cup \alpha_k.$ 
For integers $m_1,m_2,\dots,m_{|\alpha|}\in\Z_+,\,(|\alpha|= 2k)$, we
set
\[
S_\alpha(m_1,m_2,\dots,m_{|\alpha|}) = h_{a_1}(-m_1)h_{a_2}(-m_2)\cdots
h_{a_{|\alpha|}}(-m_{|\alpha|})\unit.
\]

\begin{lemma}\label{lemma:4.1.10}
Let $\alpha$ be a subset of $\{\,1,2,\dots,\ell\,\}$ with the even
cardinality
$|\alpha|$ and $|\alpha|\geq 4$. Then,
\[
S_\alpha(m_1,m_2,\dots,m_{|\alpha|})+O(\calh^+)=
(-1)^{m_1+m_2+\cdots+m_{|\alpha|}}S_\alpha+O(\calh^+).
\]
\end{lemma}
\begin{proof} We prove the lemma by
induction on
$|\alpha|$. When
$|\alpha| = 4$, it is nothing but Remark \ref{remark:4.1.9}. 
Let us suppose
$|\alpha|\geq 6$ and decompose 
$\alpha$ as $\alpha =
\alpha_1\cup\tilde{\alpha}$ where $\alpha_1 =\{a_1,a_2\,\}$ and
$\tilde\alpha=\alpha_2\cup\cdots\cup\alpha_k$ with $|\alpha_r|=2$
and $\alpha_r\cap\alpha_s=\emptyset$ for $r\ne s.$ 

Note that
\[
S_\alpha(m_1,m_2,\dots,m_{|\alpha|})=S_{a_1a_2}(m_1,m_2)*S_{\tilde{\alpha}}(
m_3,\dots,m_{|\alpha|}).
\]
By the induction hypothesis, 
\[
S_{\tilde{\alpha}}(m_3,\dots,m_{|\alpha|})\sim 
(-1)^{m_3+\cdots+m_{|\alpha|}} S_{\tilde{\alpha}}.
\]
So we have
\[
S_{\alpha}(m_1,m_2,\dots, m_{|\alpha|}) \sim
(-1)^{m_3+\cdots+m_{|\alpha|}}
(h_{a_1}(-m_1)h_{a_2}(-m_2)\unit)*S_{\alpha_2}*\cdots* S_{\alpha_k}.
\]
The proof is complete by the fact that
\[
(h_{a_1}(-m_1)h_{a_2}(-m_2)\unit)*S_{\alpha_2} =
(-1)^{m_1+m_2}S_{\alpha_1}*S_{\alpha_2},
\]
which follows from either Remark \ref{remark:4.1.9} or the induction 
hypothesis.
\end{proof}

We now use Lemma \ref{lemma:4.1.10} to prove a result similar to
Proposition \ref{proposition:4.1.6}. 

\begin{proposition}\label{lemma:4.1.11}
If $|\alpha|\geq 4$ then  as 
 a left $A[\calw]$-module or a right $A[\calw]$-module, 
 $A[\calw^{\alpha}]$ is generated by
$S_{\alpha}.$
\end{proposition}
\begin{proof}
The proof  is similar to that of Lemma \ref{lemma:4.1.5}.
For any $u\in \mathcal{W}$ and positive integers $m_1,...,m_{|\alpha|},$ set
$S_\alpha(m_1,m_2,\dots, m_{|\alpha|};u) = h_{a_1}(-m_1)h_{a_2}(-m_2)\cdots
h_{a_{|\alpha|}}(-m_{|\alpha|})u.$ Then $\calw^{\alpha}$ is
spanned by all possible $S_\alpha(m_1,m_2,\dots, m_{|\alpha|};u).$
We again use induction on $\ell(u)$ for a monomial $u$ 
to show that $S_\alpha(m_1,m_2,\dots, m_{|\alpha|};u)+O(\calh^+)$
lies in  $\Gamma_{\alpha}$ and $\Gamma_{\alpha}'$ which are the left and
right $A[\calw]$-modules generated by $S_{\alpha}+O(\calh^+),$
respectively.

If $\ell(u)=0$ then by Lemma \ref{lemma:4.1.10}, 
\[
S_\alpha(m_1,m_2,\dots, m_{|\alpha|};u)=S_\alpha(m_1,m_2,\dots, m_{|\alpha|})
\sim (-1)^{m_1+\cdots +m_{|\alpha|}}S_{\alpha}
\]
lies in $\Gamma_{\alpha}'$ and  $\Gamma_{\alpha}.$

If $\ell(u)>0$ it is clear that
\[
%\label{eqn:6}
S_\alpha(m_1,m_2,\dots,m_{|\alpha|})*u = S_\alpha(m_1,m_2,\dots,
m_{|\alpha|};u)+w
\]
where
\[
w=\sum_{n_1,\dots,n_\alpha\in \Z_+}
S_{{\alpha}}(n_1,n_2,\dots,n_{|\alpha|})
u_{n_1,...,n_{|\alpha|}},\quad
u_{n_1,...,n_{|\alpha|}}\in\calw
\]
and
$\ell(u_{n_1,...,n_{|\alpha|}})<\ell(u)$.  
Thus by the induction hypothesis, $w+O(\calh^+)$ lies in both
$\Gamma_{\alpha}'$ and  $\Gamma_{\alpha}.$

Lemma \ref{lemma:4.1.10} shows  
\[
S_\alpha(m_1,m_2,\dots,m_{|\alpha|})*u+O(\calh^+)=
(-1)^{m_1+m_2+\cdots+m_{|\alpha|}}S_{\alpha}*u+O(\calh^+)
\]
is an element of $\Gamma_{\alpha}'.$

It remains to show that $S_{\alpha}*u+O(\calh^+)\in \Gamma_{\alpha}.$ 
As in the proof of Proposition \ref{proposition:4.1.6} we have
\[
\begin{split}
S_{\alpha}*u-u*S_{\alpha} &\sim
\Res{z}\left((1+z)^{|\alpha|-1}Y(S_\alpha,z)u\right)\\
&=\sum_{i=0}^{|\alpha|-1}\binom{|\alpha|-1}{i}(S_{\alpha})_{i}u.
\end{split}
\]
Note that
\[
(S_{\alpha})_{i}=
\sum_{
\substack{m_1,....,m_{|\alpha|}\\ \sum
m_s=-|\alpha|+i+1}}
\NO h_{a_1}(m_1)h_{a_2}(m_2)\cdots
h_{a_{|\alpha|}}(m_{|\alpha|})\NO.
\]
Since $i\geq 0$ there is at least one $m_s$ positive. Thus $(S_{\alpha})_{i}u$
is a linear combination of vectors like 
$h_{a_1}(n_1)h_{a_2}(n_2)\cdots h_{a_{|\alpha|}}(n_{|\alpha|})v$ for
negative $n_s$ and a monomial $v\in \calw$ whose length is less than
the length of $u.$ By the induction hypothesis, 
$(S_{\alpha})_{i}u+O(\calh^+)\in\Gamma_{\alpha}.$ Thus 
$S_{\alpha}*u+O(\calh^+)\in \Gamma_{\alpha},$ as required. 
\end{proof}

Recall that $P$ is the collection of subsets of $\{1,...,\ell\}$ of 
even cardinalities, and that for $\alpha=\alpha_1\cup \cdots
\cup\alpha_k\in P,$ 
$S_{\alpha}=S_{\alpha_1}*\cdots *S_{\alpha_k}.$
Combining Lemma \ref{lemma:4.1.5}, Propositions 
\ref{proposition:4.1.6}, \ref{lemma:4.1.11} and 
Corollary \ref{corollary:4.1.2} we have

\begin{proposition}\label{proposition:4.1.12} The algebra 
$A(\calh^+)$ is generated by $\omega_a,$ $J_a$ 
and $h_a(-m)h_b(-n)\unit$ for $a,b =
1,\dots,\ell$ with $a\ne b$ and positive integers $m,n.$ In fact,
as a left or right $A[\calw]$-module, $A[\calw^{\alpha}]$
is generated by $h_a(-m)h_b(-n)\unit$ for
all $m,n>0$ if $\alpha=\{a,b\}$ and is generated by $S_{\alpha}$
if $|\alpha|\geq 4.$  
\end{proposition}

\begin{remark}\label{ra} In fact we can get a finite set of generators for
$A(\calh^+).$
By Proposition \ref{proposition:3.2.1} (iii), 
\[
\omega_a*S_{ab}(m,n)-S_{ab}(m,n)*\omega_a \sim mS_{ab}(m+1,n)+ mS_{ab}(m,n).
\]
for distinct $a,b$ and $m,n>0.$ Thus 
 $A[\calw^{\alpha}]$ is a generated by $S_{\alpha}(1,1)+O(\calh^+)$
as an $A[\calw]$-bimodule for $\alpha=\{a,b\}.$ In particular,
the algebra 
$A(\calh^+)$ is generated by $\omega_a$ and $J_a$ 
and $h_a(-1)h_b(-1)\unit$ for $a,b =
1,\dots,\ell$  with $a\ne b.$ 
\end{remark}

In the next two subsections we will find  a set of finite 
generators for each $A[\calw^{\alpha}]$ 
as a left or right $A[\calw]$-module for all $\alpha.$

\subsection{Consequences of the circle relation}
\label{subsection:4.2}
We derive several relations in $A(\calh^+)$ from the circle
relations. These relations 
will play important roles in the next subsection.

Recall that $S_{ab}(m,n) = h_a(-m)h_b(-n)\unit$ for distinct
$a,b$ and positive integers $m,n,$ and that $S_{ab}=S_{ab}(1,1)$.

\begin{lemma}\label{lemma:4.2.1}
For any $m,n\in\Z_+$,
\[
h_a(-1)^2S_{ab}(m,n) =
2S_{ab}(m,n)*\omega_a-2mS_{ab}(m+2,n)-2mS_{ab}(m+1,n).
\]
\end{lemma}
\begin{proof}
Proposition \ref{proposition:3.2.1} (ii) shows
\begin{align*}
S_{ab}(m,n)*\omega_a &= (L_a(-2)+L_a(-1))S_{ab}(m,n)\\
&=\frac{1}{2}h_a(-1)^2S_{ab}(m,n) +mS_{ab}(m+2,n)+mS_{ab}(m+1,n).
\end{align*}
\end{proof}

For a while, we set $u = h_a(-1)^2\unit$ and $v = h_a(-1)h_b(-n)\unit$. We
compute the result of the circle relation $u\circ v$ as follows.
Since $u_{i} = 2L_a(i-1)$, we see
\begin{align*}
u_{-2}v&= 2h_a(-2)h_a(-1)^2h_b(-n)\unit + 2h_a(-4)h_b(-n)\unit,\\
u_{-1}v&= h_a(-1)^3h_b(-n)\unit + 2h_a(-3)h_b(-n)\unit,\quad
\text{and}\\
u_{0}v&= 2h_a(-2)h_b(-n)\unit.
\end{align*}
By Proposition \ref{proposition:3.1.2} (iv),
$L(-1)(h_a(-1)^3h_b(-n)\unit)+L(0)(h_a(-1)^3h_b(-n)\unit)\sim 0$,
i.e.,
\[
3h_a(-1)^2h_a(-2)h_b(-n)\unit + n h_a(-1)^3h_b(-n-1)\unit +
(n+3)h_a(-1)^3h_b(-n)\unit\sim 0.
\]
Thus
\[
u_{-2}v\sim
-\frac{2n}{3}h_a(-1)^2S_{ab}(1,n+1)-\frac{2(n+3)}{3}h_a(-1)^2S_{ab}(1,n) +
2S_{ab}(4,n),
\]
and 
\begin{equation}\label{eqn:7}
\begin{split}
u\circ v \sim
-\frac{2n}{3}&h_a(-1)^2\bigl(S_{ab}(1,n+1)+S_{ab}(1,n)\bigr)\\
&+ 2S_{ab}(4,n) +  4S_{ab}(3,n)+2S_{ab}(2,n).
\end{split}
\end{equation}
Using Lemma \ref{lemma:4.2.1} together with (\ref{eqn:7}) yields
\[
\begin{split}
2n\bigl(S_{ab}(1,n+1)+S_{ab}(1,n)\bigr)*\omega_a
\sim &2nS_{ab}(3,n+1) +
3S_{ab}(4,n)+(2n+6)S_{ab}(3,n)\\&
+2nS_{ab}(2,n+1)+(2n+3)S_{ab}(2,n).
\end{split}
\]
Thus we have proved:

\begin{lemma}\label{lemma:4.2.2}
For distinct $a$ and $b$,
\[\begin{split}
\bigl(S_{ab}(1,m+1)+S_{ab}(1,m)\bigr)*\omega_a \sim &S_{ab}(3,m+1) +
\frac{3}{2m}S_{ab}(4,m)+\frac{m+3}{m}S_{ab}(3,m)\\&
+S_{ab}(2,m+1)+\frac{2m+3}{2n}S_{ab}(2,m).
\end{split}
\]
\end{lemma}

We also use the same argument to prove the next two lemmas. In fact,
the circle relation between $h_a(-1)h_b(-1)\unit$ and $h_a(-1)h_c(-m)\unit$
for distinct $a,b$ and $c$ with the help of relations
\begin{align*}
&L(-1)(h_a(-1)^2S_{bc}(1,m))+L(0)(h_a(-1)^2S_{bc}(1,m))\sim
0\quad\text{and}\\ &L(-1)(S_{bc}(1,m))+L(0)(S_{bc}(1,m))\sim 0
\end{align*}
gives

\begin{lemma}\label{lemma:4.2.3}
For distinct $a,b$ and $c$,
\[
\omega_a*\bigl(S_{bc}(1,m+1)+S_{bc}(1,m)\bigr)\sim
\frac{1}{2m}S_{bc}(4,m)+\frac{1}{m}S_{bc}(3,m)+\frac{1}{2m}S_{bc}(2,m).
\]
\end{lemma}

The circle relation between $h_a(-1)h_b(-1)\unit$ and
$h_a(-1)h_b(-m)\unit$ with the help of relations
\begin{align*}
&L(-1)(h_a(-1)^2S_{bb}(1,m))+L(0)(h_a(-1)^2S_{bb}(1,m))\sim 0,\\
&L(-1)(S_{bb}(1,m))+L(0)(S_{bb}(1,m))\sim 0
\end{align*}
yields

\begin{lemma}\label{lemma:4.2.4}
For distinct $a$ and $b$,
\[
\begin{split}
\omega_a*\bigl(S_{bb}(1,m+1)+S_{bb}(1,m)\bigr)&\sim
\frac{1}{2}\bigl(S_{aa}(1,m+3)+2S_{aa}(1,m+2)+S_{aa}(1,m+1)\bigr)\\
&\quad+\frac{1}{2m}\bigl(S_{bb}(4,m)+2S_{bb}(3,m)+S_{bb}(2,m)\bigr).
\end{split}
\]
\end{lemma}

%\begin{note}
%There is a difference of order of product between Lemma \ref{lemma:4.2.2}
%and Lemma \ref{lemma:4.2.3}, \ref{lemma:4.2.4}. However, we will see later
%that there are no any essential differences in the practical uses.
%\end{note}

\subsection{Finiteness}
\label{subsection:4.3}

We have already proved in 4.1 that $A[\calw^{\alpha}]$ is generated
by $S_{\alpha}$ as a left or right $A[\calw]$-module if $|\alpha|\geq 4.$
In this subsection we show that $A[\calw^{\alpha}]$ in the case $\alpha=\{a,b\}$ is generated by $h_a(-1)h_b(-m)\unit$ ($m=1,...,5$) as a left or right
$A[\calw]$-module. 

Recall that $\cals_{ab}$ is spanned by $h_a(-m)h_b(-n)\unit+O(\calh^+)$ for
positive $m,n\in \Z.$  Using the results
from the previous subsection we show that
$\cals_{ab}$ is 5 dimensional and is spanned by
$h_a(-1)h_b(-m)+O(\calh^+)$ for $m=1,...,5.$ 

Let $u = h_a(-1)h_b(-1)\unit= S_{ab}$ and $v = h_a(-1)^4\unit$.
In the following we seek the consequence of the circle relation $u\circ v$,
which turns out to be a relation in the weight $7$ space.

By direct calculations, we see
\begin{align}
u_{-2}v &= h_a(-1)^5 h_b(-2)\unit  + h_a(-2)h_a(-1)^4h_b(-1)\unit +
4h_a(-1)^3h_b(-4)\unit,\label{eqn:8}\\
u_{-1}v &= h_a(-1)^5 h_b(-1)\unit + 4h_a(-1)^3h_b(-3)\unit,\label{eqn:9}\\
u_{0}v &= 4h_a(-1)^3 h_b(-2)\unit.\label{eqn:10}
\end{align}

The important feature of the circle relation $u\circ v$ is that every term
appeared in $u\circ v$ is of the form $h_a(-1)^{2k}S_{ab}(m,n)\unit$
for $m,n\in \Z_{+}$ and $k= 1,2$. 
The case $k=1$ was already considered in Lemma \ref{lemma:4.2.1}
and the following relation was obtained:
\begin{equation}\label{eqn:11}
h_a(-1)^2S_{ab}(m,n) =
2S_{ab}(m,n)*\omega_a-2mS_{ab}(m+2,n)-2mS_{ab}(m+1,n).
\end{equation}

Now we turn to the case $k=2$.

\begin{lemma}\label{lemma:4.3.1}
For distinct $a$ and $b$,
\begin{align*}
&h_a(-1)^4S_{ab}(1,m)\\
&\sim 4S_{ab}(1,m)*\omega_a^2\\
&\quad -\bigl(16S_{ab}(3,m) +
4S_{ab}(2,m)-4mS_{ab}(1,m+1)-4(m+3)S_{ab}(1,m)\bigr)*\omega_a\\
&\quad +36S_{ab}(5,m)+36S_{ab}(4,m)-4mS_{ab}(3,m+1)\\
&\quad -4m S_{ab}(2,m+1)-4(m+3)S_{ab}(3,m)-4(m+3)S_{ab}(2,m).
\end{align*}
\end{lemma}
\begin{proof}
Set $w = h_a(-1)^2S_{ab}(1,m)$. Since
\begin{align*}
w*\omega_a&\sim (L_a(-2)+L_a(-1))w\\
&=\frac{1}{2}h_a(-1)^4S_{ab}(1,m) + 3h_a(-3)h_a(-1)^2h_b(-m)\unit +
L_a(-1)w,
\end{align*}
we have
\begin{equation}\label{eqn:12}
h_a(-1)^4S_{ab}(1,m)\sim 2 w*\omega_a - 6
h_a(-3)h_a(-1)^2h_b(-m)\unit - 2L_a(-1)w.
\end{equation}
Note that $L(-1)w=(L_a(-1)+L_b(-1))w \sim -L(0)w$, i.e.,
\begin{equation}\label{eqn:13}
L_a(-1)w\sim -mh_a(-1)^2S_{ab}(1,m+1) - (m+3)h_a(-1)^2S_{ab}(1,m).
\end{equation}
Combining (\ref{eqn:12}) and
(\ref{eqn:13}) together proves 
\[
\begin{split}
h_a(-1)^4S_{ab}(1,m)
&\sim 2h_a(-1)^2S_{ab}(1,m)*\omega_a - 6h_a(-1)^2S_{ab}(3,m)\\
&\quad +2mh_a(-1)^2S_{ab}(1,m+1)+2(m+3)h_a(-1)^2S_{ab}(1,m).
\end{split}
\]
Finally, using (\ref{eqn:11}) we get the desired result.
\end{proof}

Note that each term in (\ref{eqn:8}), (\ref{eqn:9}) and (\ref{eqn:10})
except $h_a(-2)h_a(-1)^4h_b(-1)\unit$
has the form $h_a(-1)^{2k}S_{ab}(1,m)$ for $k = 1,2$.
Using the relation $L(-1)(h_a(-1)^4S_{ab}(1,1)) +
L(0)(h_a(-1)^4S_{ab}(1,1))\sim 0$,
i.e.,
\[
h_a(-1)^4S_{ab}(2,1)\sim -\frac{1}{5}h_a(-1)^4S_{ab}(1,2)-\frac{6}{5}
h_a(-1)^4S_{ab}(1,1),
\]
we see that $h_a(-1)^4S_{ab}(2,1)$ is a linear
combination of the elements $S_{ab}(m,n)*\omega_a^i.$
Thus we have by
(\ref{eqn:11}) and Lemma
\ref{lemma:4.3.1}:

\begin{lemma}\label{lemma:4.3.2} 
Let $u = S_{ab}(1,1)$ and $v = h_a(-1)^4\unit$. Then $u_{-j}v+O(\calh^+)$
for $j=0,1,2$ is a linear combination of $S_{ab}(m,n)*\omega_a^i+O(\calh^+)$ 
where $0\leq i\leq 2$ and $m+n\leq 5+j-2i$. In particular,
$u\circ v$ is a
linear combination of the elements $S_{ab}(m,n)*\omega_a^i$ modulo 
$O(\calh^+)$ where $0\leq i\leq 2$ and
$m+n\leq 7-2i$. 
\end{lemma}

For a positive integer $k$, let $\mathcal{S}_{ab}(k)$ be the span 
of $S_{ab}(m,n)+O(\calh^+)$ for $m,n\in\Z_{+}$
such that $m+n\leq k+1$. 

\begin{lemma}\label{lemma:4.3.3}
The $\cals_{ab}(k)$ is spanned by 
$\{S_{ab}(1,m)+O(\calh^+)\,|\,1\leq m\leq k\}$.
In particular, the dimension of $\cals_{ab}(k)$ is less than or equal to 
$k.$ 
\end{lemma}
\begin{proof}
We prove the lemma by induction on $k$. When $k=1$, it is
clear. Let $k\geq 2$ and $m,n>0$ with $m+n=k+1.$  
Then the relation $L(-1)S_{ab}(m,n)+L(0)S_{ab}(m,n)\sim 0$ gives
\[
mS_{ab}(m+1,n)+nS_{ab}(m,n+1)+(m+n)S_{ab}(m,n)\sim 0,
\]
which implies 
\[
mS_{ab}(m+1,n)+nS_{ab}(m,n+1)\equiv 0\quad \mod \quad \cals_{ab}(k-1).
\]
Thus
\begin{equation}\label{eqn:14}
S_{ab}(n,k-n+1) \equiv (-1)^{n-1}\binom{k-1}{n-1}S_{ab}(1,k)\quad\mod\quad
\cals_{ab}(k-1).
\end{equation}
This shows that $\dim \cals_{ab}(k)/\cals_{ab}(k-1)\leq 1.$
The induction hypothesis then yields that 
$\cals_{ab}(k)$ is spanned by 
$\{S_{ab}(1,m)+O(\calh^+)\,|\,1\leq m\leq k\}.$ 
\end{proof}

We go back to the circle relation $u\circ v$ for $u = S_{ab}(1,1)$ and $v =
h_a(-1)^4\unit.$ By Lemmas \ref{lemma:4.3.2} and \ref{lemma:4.3.3}, we can
write $u\circ v$ as
\[
u\circ v \sim \sum_{\substack{j=0,1,2\\
m\leq 6-2j}}x_{mj}S_{ab}(1,m)*\omega_a^j
\]
with some scalars $x_{mj}$. Using Lemma \ref{lemma:4.2.2} and Lemma
\ref{lemma:4.3.3} proves
\begin{equation}\label{eqn:15}
u\circ v \sim \sum_{i=1,2}x_iS_{ab}(1,1)*\omega_a^i +
\sum_{m=1}^6y_mS_{ab}(1,m)
\end{equation}
with some scalars $x_i$ and $y_m$. We call this expression of $u\circ
v$ the \textit{normal form}. The same process also shows that each 
term (a homogeneous vector) 
occurring in (\ref{eqn:8})-(\ref{eqn:10}) has a normal form 
like (\ref{eqn:15}) and the weights of homogeneous vectors in
the normal form are less than or equal to the weight of the
original term. 

Since $u\circ v$  acts trivially on each top level
of $\calh^+$-modules, the right hand side of
(\ref{eqn:15}) acts also trivially on each top level
of $\calh^+$-modules. Let $\l=\sum_{i=1}^{\ell}\l_ih_i\in\frakh.$ 
Then the top level of $\calh^+$-module $M(1,\lambda)$
is one dimensional.  Now $S_{ab}(1,m)$ acts on the top level of 
$M(1,\lambda)$ as $S_{ab}(1,m)= (-1)^{m+1}\lambda_a\lambda_b$ and
$\omega_a$ acts on this space as $\omega_a = \frac{1}{2}\lambda_a^2$:
See Table 1 at the end of this subsection. Then the
evaluation of the right hand side of (\ref{eqn:15}) 
on the top level of this module shows
\[
\sum_{i=1,2}2^{-i}x_i\lambda_a^{2i+1}\lambda_b  +
\sum_{m=1}^6y_m(-1)^m\lambda_a\lambda_b = 0.
\]
Since $\l$ is arbitrary we immediately see that $x_i = 0$ for $i=1,2$  
and obtain the relation $\sum_{m=1}^6y_mS_{ab}(1,m)\sim 0$.

Now, suppose that there exists $m$ such that $y_m\neq 0$. Then we have a
(nontrivial) relation  among $S_{ab}(1,m),\,(1\leq m\leq 6)$. The
coefficients $y_m$ might be zero for all $m.$ However, we are able to show that
$y_6\neq 0.$

{}From Lemma \ref{lemma:4.3.2} and the proof of Lemma \ref{lemma:4.3.3},
there is no contribution to $S_{ab}(1,6)$  in the normal form from
either
$u_{-1}v$ or $u_0v.$ To find $y_6$, it is enough to consider
the term $u_{-2}v$ in $u\circ v.$ We define an equivalence 
relation  $\succsim $ on $\calh^+$ such that $w_1\succsim w_2$
if and only if there exists $v$ whose homogeneous
components have weight less than $7$ 
such that $w_1-w_2\sim v.$  

Let us recall that
\begin{align*}
&u_{-2}v = h_a(-1)^5h_b(-2)\unit  + h_a(-2)h_a(-1)^4h_b(-1)\unit +
4h_a(-1)^3h_b(-4)\unit,\quad\text{and}\\
&h_a(-1)^4S_{ab}(2,1)\sim -\frac{1}{5}h_a(-1)^4S_{ab}(1,2)-\frac{6}{5}
h_a(-1)^4S_{ab}(1,1).
\end{align*}
Then
\begin{equation}\label{eqn:16}
u_{-2}v\succsim \frac{4}{5}h_a(-1)^4S_{ab}(1,2) + 4h_a(-1)^2S_{ab}(1,4).
\end{equation}

{}From Lemma \ref{lemma:4.3.1} we see that 
\begin{equation}\label{eqn:17}
h_a(-1)^4S_{ab}(1,2)\succsim 4S_{ab}(1,2)*\omega_a^2 - 16
S_{ab}(3,2)*\omega_a + 36 S_{ab}(5,2).
\end{equation}
By (\ref{eqn:11}) we know that
\begin{equation}\label{eqn:18}
h_a(-1)^2S_{ab}(1,4)\succsim 2S_{ab}(1,4)*\omega_a - 2S_{ab}(3,4).
\end{equation}
Using (\ref{eqn:16}), (\ref{eqn:17}) and (\ref{eqn:18}) shows
\begin{equation}
\begin{split}\label{eqn:19}
u_{-2}v&\succsim \frac{16}{5}S_{ab}(1,2)*\omega_a^2\\
&\quad +\bigl(8S_{ab}(1,4)-\frac{64}{5}S_{ab}(3,2)\bigr)*\omega_a
+\frac{144}{5}S_{ab}(5,2) - 8S_{ab}(3,4).
\end{split}
\end{equation}

Using  the relation
$S_{ab}(1,2)*\omega_a^2\succsim 
\bigl(S_{ab}(3,2)+\frac{3}{2}S_{ab}(4,1)\bigr)*\omega_a$ which follows from 
Lemma \ref{lemma:4.2.2}, we reduce (\ref{eqn:19}) to 
\[
%\label{eqn:20}
\begin{split}
u_{-2}v&\succsim \bigl(
-\frac{48}{5}S_{ab}(3,2)
+\frac{24}{5}S_{ab}(4,1)+8S_{ab}(1,4)\bigr)*\omega_a\\
&\quad+\frac{144}{5}S_{ab}(5,2)-8S_{ab}(3,4)
\end{split}
\]
which, by (\ref{eqn:14}), is equivalent to 
\[
u_{-2}v\succsim  -\frac{128}{5}S_{ab}(1,4)*\omega_a
+\frac{144}{5}S_{ab}(5,2) - 8S_{ab}(3,4).
\]
To find $y_6$, we still need to reduce the term $S_{ab}(1,4)$ into the
normal form. Again by using Lemma \ref{lemma:4.2.2} we see that
\[
S_{ab}(1,4)*\omega_a\succsim S_{ab}(3,4)+\frac{1}{2}S_{ab}(4,3).
\]
Finally we have $u_{-2}v\succsim -64S_{ab}(1,6)$ 
by (\ref{eqn:14}). This shows that $y_6 = -64$.

Thus we have proved:

\begin{lemma}\label{lemma:4.3.4}
In $A(\calh^+)$, the element $S_{ab}(1,6)+O(\calh^+)$ is a linear
combination of the
elements $S_{ab}(1,m)+O(\calh^+),\ (1\leq m\leq 5).$
\end{lemma}

Now we prove one of the important results in this paper. Recall
that $\cals_{ab}$ is spanned by   all $S_{ab}(m,n)$, $(m,n\in\Z_{+})$.

\begin{proposition}\label{proposition:4.3.5}
We have  $\dim \cals_{ab}=5$ with a basis $S_{ab}(1,m)+O(\calh^+),
\,(1\leq m\leq 5).$ 
\end{proposition}
\begin{proof}
We first prove $\dim A(\cals_{ab})\leq 5$. It is enough to
show that $S_{ab}(m,n)+O(\calh^+)$ for $m+n\geq 7$ is expressed as a linear
combination of $S_{ab}(1,m)+O(\calh^+),\,(1\leq m\leq 5)$. Note from 
 Proposition \ref{proposition:3.2.1} (iii) that
\begin{equation}\label{eqn:21}
mS_{ab}(1,m+1) \sim [\omega_b,S_{ab}(1,m)]-mS_{ab}(1,m)
\end{equation}
for any positive integer $m$ 
where $[\omega_b,S_{ab}(1,m)]=\omega_b*S_{ab}(1,m)-S_{ab}(1,m)*\omega_b.$
Since $S_{ab}(1,6)+O(\calh^+)$ is a linear combination of
$S_{ab}(1,m)+O(\calh^+),\,(1\leq m\leq 5)$ by Lemma \ref{lemma:4.3.4}, 
the element $[\omega_b,S_{ab}(1,5)]$ is also a linear combination of 
$S_{ab}(1,m)+O(\calh^+)$ for $m=1,\dots,5.$ 
We also see from (\ref{eqn:21}) that the space spanned by
$S_{ab}(1,m)+O(\calh^+),\,(1\leq m\leq 5)$ is closed under the bracket
$[\omega_b+O(\calh^+),\,\cdot\,].$ 

We now prove by induction on  positive integer $n$ greater than or equal to
$6$ that $S_{ab}(1,n)+O(\calh^+)$ is a linear combination of 
$S_{ab}(1,m)+O(\calh^+),\,(1\leq m\leq 5).$ Suppose that 
this is true for $6,\dots, n$. Then since $S_{ab}(1,n)+O(\calh^+)$ is a
linear combination of $S_{ab}(1,m)+O(\calh^+),\,(1\leq m\leq 5).$ From the 
discussion in the previous paragraph, $[\omega_b,S_{ab}(1,n)]+O(\calh^+)$ 
is also a linear combination of
$S_{ab}(1,m)+O(\calh^+),\,(1\leq m\leq 5)$ and so is 
$S_{ab}(1,n+1)+O(\calh^+).$ 

It remains to show that $S_{ab}(1,m)+O(\calh^+),\,(1\leq m\leq 5)$ are linearly
independent. Suppose that there exist scalars $x_m,\,(1\leq m\leq 5)$ such
that 
\begin{equation}\label{eqn:22}
\sum_{m=1}^5 x_mS_{ab}(1,m)+O(\calh^+) = 0.
\end{equation}
We have to prove that $x_1= x_2=\dots=x_5=0.$ This is achieved by
the following evaluation method. We
evaluate (\ref{eqn:22}) on the top levels of irreducible
modules $\calh^-$, $M(1,\lambda)$ and $\calh(\theta)^-$.
Consider the basis $\{h_a(-1)\unit|a=1,\dots,\ell\}$ for
$\calh^-$ and the basis $\{h_a(-1/2)\unit|a=1,\dots,\ell\}$
for $\calh(\theta)^-$.  Let $E_{ab}$ be the corresponding matrix 
element with respect to the basis. 
Also write $\l=\sum_{i=1}^{\ell}\l_ih_i$.
Then the action of $S_{ab}(1,m),\,(1\leq m\leq 5)$ on the top levels
of $\calh^-$, $M(1,\lambda)$ and $\calh(\theta)^-$ are listed in Table
1; $S_{ab}(1,m)$ acts trivially on the top levels of the modules
$\calh^+$ and $\calh(\theta)^+$.

\vskip 5ex

\begin{center}
\begin{tabular}{|c|c|c|c|}
\hline
{}&$\calh^-$&$M(1,\lambda)$&$\calh(\theta)^-$\\
\hline
$S_{ab}(1,1)$&$E_{ab}+E_{ba}$&$\lambda_a\lambda_b$&$\frac{1}{2}(E_{ab}+E_{ba
})$\\
\hline
$S_{ab}(1,2)$&$-2E_{ab}$&$-\lambda_a\lambda_b$&$-\frac{3}{4}E_{ab}-\frac{1}{
4}E_{ba}$\\
\hline
$S_{ab}(1,3)$&$3E_{ab}$&$\lambda_a\lambda_b$&$\frac{15}{16}E_{ab}+\frac{3}{1
6}E_{ba}$\\
\hline
$S_{ab}(1,4)$&$-4E_{ab}$&$-\lambda_a\lambda_b$&$-\frac{35}{32}E_{ab}-\frac{5
}{32}E_{ba}$\\
\hline
$S_{ab}(1,5)$&$5E_{ab}$&$\lambda_a\lambda_b$&$\frac{315}{256}E_{ab}+\frac{35
}{256}E_{ba}$\\
\hline
\end{tabular}

\vskip 5ex
Table 1
\end{center}

Then we have the linear system
\begin{align*}
&\phantom{-}x_1 = 0,\\
&-2x_2+3x_3-4x_4+5x_5 = 0,\\
&-x_2+x_3-x_4+x_5=0,\\
&-\frac{3}{4}x_2+\frac{15}{16}x_3-\frac{35}{32}x_4+\frac{315}{256}x_5 = 0,\\
&-\frac{1}{4}x_2+\frac{3}{16}x_3-\frac{5}{32}x_4+\frac{35}{256}x_5 = 0.
\end{align*}
The only solution to this linear system 
is $x_1=\cdots =x_5=0,$ as expected.
\end{proof}

Combining Propositions \ref{proposition:4.1.12} and 
\ref{proposition:4.3.5} gives
\begin{corollary} Let $\alpha=\{a,b\},\,(a\neq
b)$. Then as a left or right $A[\calw]$-module, 
$A[\calw^{\alpha}]$ is generated by  
$S_{ab}(1,m),\,(1\leq m\leq 5).$ 
\end{corollary}

The importance of the elements of $S_{ab}(1,m),\,(1\leq m\leq 5)$
will become clear in the next section.

\section{Matrix subalgebras}
\label{section:5}
In this section we use the elements $S_{ab}(1,m)+O(\calh^+),\,(1\leq m\leq 5)$
to construct two ideals in $A(\calh^+).$ Both 
are isomorphic to the $\ell\times \ell$ matrix algebra. 
These two ideals characterize the modules $\calh^-$ and $\calh(\theta)^-$.

\subsection{Matrix elements}
\label{subsection:5.1}
Motivated by Table 1, we define elements
$E_{ab}^u,\bar{E}_{ba}^u, E_{ab}^t,\bar{E}_{ba}^t$ and $\Lambda_{ab}$ as
follows; for distinct $a$ and $b$, set 

\begin{align}
&E_{ab}^u = 5S_{ab}(1,2)+25S_{ab}(1,3)+36S_{ab}(1,4)+16S_{ab}(1,5),\notag\\
&\bar{E}_{ba}^u =
S_{ab}(1,1)+14S_{ab}(1,2)+41S_{ab}(1,3)+44S_{ab}(1,4)+16S_{ab}(1,5),\notag\\
&E_{ab}^t = -16\bigl(
3S_{ab}(1,2)+14S_{ab}(1,3)+19S_{ab}(1,4)+8S_{ab}(1,5)\bigr),\label{eqn:23}\\
&\bar{E}_{ba}^t = -16\bigl(
5S_{ab}(1,2)+18S_{ab}(1,3)+21S_{ab}(1,4)+8S_{ab}(1,5)\bigr),\notag\\
&\Lambda_{ab}=45S_{ab}(1,2)+190S_{ab}(1,3)+240S_{ab}(1,4)+96S_{ab}(1,5).\notag
\end{align}

The action of these elements on the top levels of 
of irreducible modules $\calh^-,$ $M(1,\lambda)$ and $\calh(\theta)^-$ is
listed in Table 2.

\vskip 5ex

\begin{center}
\begin{tabular}{|c|c|c|c|}
\hline
{}&$\calh^-$&$M(1,\lambda)$&$\calh(\theta)^-$\\
\hline
$E_{ab}^u $&$E_{ab}$&$0$&$0$\\
\hline
$\bar{E}_{ba}^u$&$E_{ba}$&$0$&$0$\\
\hline
$E_{ab}^t$&$0$&$0$&$E_{ab}$\\
\hline
$\bar{E}_{ba}^t$&$0$&$0$&$E_{ba}$\\
\hline
$\Lambda_{ab}$&$0$&$\lambda_a\lambda_b$&$0$\\
\hline
\end{tabular}
\vskip 5ex
Table 2
\end{center}

\vskip 5ex

Recall that $[u]=u+O(\calh^+).$ 
\begin{lemma}\label{lemma:5.1.1}
For distinct $a$ and $b$,
\[
[E_{ba}^u]= [\bar{E}_{ba}^u],\quad [E_{ba}^t] =
[\bar{E}_{ba}^t]\quad\text{and}\quad [\Lambda_{ab}] = [\Lambda_{ba}].
\]
\end{lemma}
\begin{proof}
One could  prove the lemma using the fact that $h_a(-n)h_b(-1)\unit$ is a
linear combination of the elements $S_{ab}(1,m)$'s with 
a lengthy computation. Here we give a short proof by evaluation method 
discussed in the proof of Proposition \ref{proposition:4.3.5}.

Consider $D_{ab} = E_{ba}^u - \bar{E}_{ba}^u$ or 
$E_{ba}^t -\bar{E}_{ba}^t$ or $\Lambda_{ab}-\Lambda_{ba}$.
Clearly $D_{ab}$ is a linear combination of $S_{ab}(1,m)$ for
$m=1,\dots,5$ modulo $O(\calh^+)$ 
and acts trivially on top levels of irreducible modules in Table 2. 
The proof of (\ref{eqn:22}) shows that $D_{ab} = 0$. 
\end{proof}

Thanks to Lemma \ref{lemma:5.1.1}, we use $E_{ba}^u$ and $E_{ba}^t$ instead
of $\bar{E}_{ba}^u$ and $\bar{E}_{ba}^t$. 

\begin{remark}\label{remark:5.1.2}
Solving the linear system (\ref{eqn:23}) for $S_{ab}(1,m)$'s gives rise to
\begin{align}
S_{ab}(1,1)&= E_{ab}^u+E_{ba}^u + \Lambda_{ab} + \frac{1}{2}E_{ab}^t +
\frac{1}{2}E_{ba}^t,\notag\\
S_{ab}(1,2)&= -2 E_{ab}^u- \Lambda_{ab} - \frac{3}{4}E_{ab}^t -
\frac{1}{4}E_{ba}^t,\notag\\
S_{ab}(1,3)&= 3 E_{ab}^u + \Lambda_{ab} + \frac{15}{16}E_{ab}^t +
\frac{3}{16}E_{ba}^t,\label{eqn:n24}\\
S_{ab}(1,4)&= -4 E_{ab}^u- \Lambda_{ab} - \frac{35}{32}E_{ab}^t -
\frac{5}{32}E_{ba}^t,\notag\\
S_{ab}(1,5)&= 5 E_{ab}^u + \Lambda_{ab} + \frac{315}{256}E_{ab}^t +
\frac{35}{256}E_{ba}^t.\notag
\end{align}
This is essentially Table 1.

\end{remark}

\begin{remark}\label{remark:5.1.3} By Remarks \ref{ra} and 
\ref{remark:5.1.2}, $A(\calh^+)$ is generated by $\omega_a,J_a,
E_{ab}^u, E_{ab}^t$ and $\Lambda_{ab}$ for $a,b=1,...,\ell$ 
and $a\ne b.$ 
\end{remark}

\subsection{Relations among $\omega_a, E_{ab}^{u}$ and $E_{ab}^{t}$}
\label{subsection:5.2}
In this subsection we  use the evaluation method developed in the preceding
subsections to derive some relations among $\omega_a, E_{ab}^{u}$
and $E_{ab}^{t}.$

First we give the actions of $\omega_a$ and $J_a$ on top levels 
of the known irreducible modules for $\calh^+$ in Table 3. Here we use $I$ to
denote the $\ell\times\ell$ identity matrix. 

\vskip 5ex
\begin{center}
\begin{tabular}{|c|c|c|c|c|c|}
\hline
&$\calh^+$&$\calh^-$&$M(1,\lambda)$&$\calh(\theta)^+$&$\calh(\theta)^-$\\
\hline
$\omega_a$&$0$&$E_{aa}$&$\frac{1}{2}\lambda_a^2$&$\frac{1}{16}$&$\frac{1}{16
}I+\frac{1}{2}E_{aa}$\\
\hline
$J_a$&$0$&$-6E_{aa}$&$\lambda_a^4-\frac{1}{2}\lambda_a^2$&$\frac{3}{128}$&
$\frac{3}{128}I-\frac{3}{8}E_{aa}$\\
\hline
\end{tabular}
\end{center}
\vskip 5ex
\begin{center}
Table 3
\end{center}

For convenient we shall use  $E_{ab}^*$ for either
$E_{ab}^u$ or $E_{ab}^t$. Whenever $E_{ab}^*$ appears in equalities, it is
supposed to be the same element.

Recall that ${\cal S}_{ab}$ is the linear span of
$h_a(-m)h_b(-n)\unit+O(\calh^+)$
for $m,n\in\Z_+$. Then ${\cal S}_{ab}$ is 5-dimensional with 
a basis $\{[S_{ab}(1,m)]+O(\calh^+)|m=1,\dots,5\}$. 
Let $\calm_{ab}^u$ be the vector subspace of $A(\calh^+)$ spanned by 
the elements $[E_{ab}^u], [E_{ba}^u]$ and $\calm_{ab}^t$ be
the vector subspace spanned by $[E_{ab}^t], [E_{ba}^t]$, and $\calm_{ab}
=
\calm_{ab}^u\oplus\calm_{ab}^t$.

The following lemma is frequently used in this subsection.

\begin{lemma}\label{lemma:5.2.1a}
{\rm (i)} If an element $E \in \cals_{ab}$ acts trivially on the top
levels of $\calh^\pm,$ $M(1,\lambda),$ $\lambda(\neq 0)\in\frakh$ and
$\calh(\theta)^\pm$, then $E= 0$.

\noindent
{\rm (ii)} If an element $E\in\calm_{ab}$ acts trivially on the top
levels of the modules $\calh^-$ and $\calh(\theta)^-$, then $E= 0$.
\end{lemma}
\begin{proof}
Write $E$ as a linear combination of the elements $E_{ab}^u+O(\calh^+),
E_{ba}^u+O(\calh^+),
E_{ab}^t+O(\calh^+), E_{ba}^t+O(\calh^+)$ and
$\Lambda_{ab}+O(\calh^+)$. Then the result is immediate.
\end{proof}

Now we prove relations between $\omega_a$ and $E_{bc}^*$ in $A(\calh^+).$
As we pointed out before, we sometimes identify $u\in \calh^+$ with
its image $[u]=u+O(\calh^+)$ in $A(\calh^+)$ when the content is clear.
For short we sometimes write $xy$ for $x*y$ in $A(\calh^+).$ 

\begin{lemma} \label{lemma:5.2.2a}
For any $a$ and distinct $b,c$, 
\begin{align*}
&\omega_a* E_{bc}^u = \delta_{ab}E_{bc}^u,\quad E_{bc}^u*\omega_a
=\delta_{ac}E_{bc}^u
\quad\text{and}\\
&\omega_a*E_{bc}^t = (\frac{1}{16}+\frac{1}{2}\delta_{ab})E_{bc}^t,\quad
E_{bc}^t *\omega_a= (\frac{1}{16}+\frac{1}{2}\delta_{ac})E_{bc}^t
\end{align*}
in $A(\calh^+).$ 
\end{lemma}
\begin{proof}
Since all equalities in the lemma hold on any top levels of the modules
$\calh^-$ and  $\calh(\theta)^-$, it suffices to prove by Lemma 
\ref{lemma:5.2.1a} that the left hand 
sides of the equalities belong to $\calm_{ab}$.

First we consider the case that
$a,b$ and $c$ are all distinct.  Then by Lemma \ref{lemma:4.2.3},
Proposition \ref{proposition:4.3.5} and Remark \ref{remark:5.1.2}, we
know that
$\omega_a* E_{bc}^*$ is expressed as a linear combination of the elements
$\omega_a*S_{bc}(1,1)$ and
$E_{bc}^u, E_{cb}^u,E_{bc}^t, E_{cb}^t$ and $\Lambda_{bc}$ modulo
$O(\calh^+):$ 
\begin{equation}\label{eqn:5.2.1a}
\omega_a* E_{bc}^* = x_0\omega_a*S_{bc}(1,1) + x_1 E_{bc}^u+x_2
E_{cb}^u+x_3E_{bc}^t+x_4E_{cb}^t+x_5\Lambda_{bc}.
\end{equation}
for some scalars $x_i.$ 
Evaluation of (\ref{eqn:5.2.1a}) on the top level of $M(1,\lambda)$ shows
\[
\frac{1}{2}x_0\lambda_a^2\lambda_b\lambda_c + x_5\lambda_b\lambda_c = 0
\quad\text{for all}\quad 0\neq \lambda=\sum_{i}\l_ih_i\in\frakh,
\]
which implies $x_0 =x_5=0$. Therefore $\omega_a*E_{bc}^*$ is an element
of $\calm_{bc}$. Further, since $[\omega_a, E_{bc}^* ]\sim
(L_a(0)+L_a(-1))E_{bc}^*=0$, we have
$E_{bc}^**\omega_a\in\calm_{bc}$.

Next we consider the case either $a=b$ or $a=c$.
 Lemma \ref{lemma:4.2.2} and Remark \ref{remark:5.1.2} show that
$E_{bc}^**\omega_a$ is a linear combination of $S_{bc}(1,1)*\omega_a,
E_{bc}^u, E_{cb}^u,E_{bc}^t, E_{cb}^t$ and $\Lambda_{bc}$ modulo $O(\calh^+):$
\[
%\label{eqn:25}
E_{bc}^**\omega_a = x_0 S_{bc}(1,1)*\omega_a + x_1 E_{bc}^u+x_2
E_{cb}^u+x_3E_{bc}^t+x_4E_{cb}^t+x_5\Lambda_{bc}
\]
with some scalars $x_i,\,(0\leq i\leq 5)$. The same evaluation 
in the previous paragraph shows  $x_0 = x_5 = 0$, i.e.,
$E_{bc}^*\omega_a\in \calm_{ab}$. Note that $\omega_a*E_{bc}^* =
E_{bc}^**\omega_a + [\omega_a,E_{bc}^*]$ and $[\omega_a,E_{bc}^*] \sim
L_a(-1)E_{bc}^*+ L_a(0)E_{bc}^*$. Since the space
spanned by $S_{bc}(m,n)$ for $m,n>0$ is invariant under 
$L_a(-n)$  for $n=0,1,$ $\omega_a*E_{bc}^*+O(\calh^+)$
lies in ${\cal S}_{bc}.$ The same evaluation then concludes
that $\omega_a*E_{bc}^*+O(\calh^+)\in
\calm_{bc}$.
\end{proof}

Now we turn our attention to $J_a.$ We have
\begin{proposition}\label{proposition:5.2.3} For any $a$ and distinct
$b,c$, the spaces
$\calm_{bc}^u$ and $\calm_{bc}^t$ are closed under the left and right
multiplication by $J_a$.
\end{proposition}
\begin{proof} We only consider the case
that $a=b$ or $a=c$ here. The case that $a,b$ and $c$ are all distinct is
treated in Subsection \ref{subsection:5.3}; see Corollary 
\ref{corollary:5.3.5a}. 

We first prove that if $a=b$ then $J_a*S_{bc}(1,m),$ is
 a linear combination of
$S_{ac}(1,n)*\omega_a^j,\,(j\in\N,1\leq n\leq 5)$ modulo
$O(\calh^+).$ By Proposition \ref{proposition:4.3.5}, it is enough
to show that $J_a*S_{ac}(1,m),$ is
a linear combination of
$S_{ac}(p,q)*\omega_a^j,\,(j\geq 0, p,q>0)$ modulo
$O(\calh^+).$

Recall that  
\begin{align*}
J_a*S_{ac}(1,m) &= \sum_{i=0}^4\binom{4}{i}(J_a)_{i-1}S_{ac}(1,m),\quad (1\leq
m\leq 5)\\
&=\sum_{i=0}^4\binom{4}{i}h_c(-m)(J_a)_{i-1}h_a(-1)\unit,
\end{align*}
where we set $J_a(z) = Y(J,z) = \sum_{n\in\Z}(J_a)_{n}z^{-n-1}$.
The space $\calh_a^-$ is a Virasoro module with respect to a Virasoro
element $\omega_a$ and decomposes into the direct sum of irreducible
modules $\calh_a^- = \oplus_{i=0}^\infty L(1,(2i+1)^2)$ 
where $L(c,h)$
denotes an irreducible highest weight module with highest
weight $h$ and central charge $c$ for the Virasoro element $\omega_a$ 
(cf.~\cite{DG}). Since 
$\wt{(J_a)_{i-1}h_a(-1)\unit}\leq 5,$
 we see $(J_a)_{i-1}h_a(-1)\unit\in
L(1,1)$ for $0\leq i\leq 4$.  Therefore, $J_a*S_{ac}(1,m)$ is
a linear combination of the elements
\[
%\label{eqn:26}
L_a(-n_1)\cdots L_a(-n_k)S_{ac}(1,m),\quad n_i\geq 1.
\]

It is enough to show that $L_a(-n_1)\cdots L_a(-n_k)S_{ac}(s,t),
\,(n_i\geq 1, s,t>0)$ are linear
combinations of the elements $S_{ac}(p,q)*\omega_a^j,\,
(j\geq 0, p,q>0)$ modulo $O(\calh^+).$  
By Proposition \ref{proposition:3.2.1} (i) we can
assume that $n_1\leq 2$. Set $v = L_a(-n_2)\cdots
L_a(-n_k)S_{ac}(s,t).$ 
Then  by Proposition \ref{proposition:3.2.1} (ii),
$L_a(-2)v \sim v*\omega_a - L_a(-1)v.$  Note that
\begin{align*}
L_a(-1)v &=\sum_{i}(n_i-1)L_a(-n_2)\cdots L_a(-n_i-1)\cdots
L_a(-n_k)S_{ac}(s,t)\\ &\ \ \ \ +sL_a(-n_2)\cdots
L_a(-n_k)S_{ac}(s+1,t).
\end{align*}
Induction on $k$ gives the desired result.

Since 
\[
J_a*S_{ac}(1,m)-S_{ac}(1,m)*J_a\sim\sum_{i=0}^3\binom{3}{i}(J_a)_iS_{ac}(1,m
),\quad
(1\leq m\leq 5),
\]
the same argument above shows that $J_a*S_{ac}(1,m)-S_{ac}(1,m)*J_a$
and thus $S_{ac}(1,m)*J_a$ are also linear combinations of
$S_{ac}(1,n)*\omega_a^j,\,(j\in\N,1\leq n\leq 5)$ modulo
$O(\calh^+).$

If $a=c,$  by Proposition \ref{proposition:4.3.5}, the span 
of $S_{ba}(1,m)+O(\calh^+)$ for $m=1,...,5$ is the same as the span 
$S_{ba}(n,1)+O(\calh^+)$ for $n=1,...,5.$ This implies
that both $J_a*S_{ba}(1,m)$ and $S_{ba}(1,m)*J_a$ are 
linear combinations of
$S_{ba}(1,n)*\omega_a^j,\,(j\in\N,1\leq n\leq 5)$ modulo
$O(\calh^+).$

Then  by Lemma \ref{lemma:5.2.2a},
$J_a*S_{bc}(1,m)$ and $S_{bc}(1,m)*J_a$ are linear combinations of
$E_{bc}^u, E_{cb}^u,E_{bc}^t, E_{cb}^t$ and $\Lambda_{bc}*\omega_a^j$
modulo $O(\calh^+),$ and so are
$J_a*E_{bc}^*$ and $E_{bc}^**J_a$. Now, let
\[
J_a*E_{bc}^* = \sum_{j\in\N}x_{j}\Lambda_{bc}\omega_a^j +
y_1E_{bc}^u+y_2E_{cb}^u+y_3E_{bc}^t+y_4E_{cb}^t
\]
modulo $O(\calh^+)$ with some scalars $x_{j}$ and $y_i$.
Then the evaluation of these expressions on the top level of 
$M(1,\lambda)$ shows $x_{j} = 0$ and then $J_a*E_{bc}^*\in\calm_{bc}^*$.
Similarly we can prove $E_{bc}^**J_a\in\calm_{bc}^*$. Namely, $\calm_{bc}^*$
is closed under both right and left multiplications by $J_a$. 
\end{proof}

\subsection{Relations among $E_{ab}^*$ }
\label{subsection:5.3}

For convenience, we set
$\bar{\calm}_{ab} =\cals_{ab} =
\calm_{ab}\oplus\C(\Lambda_{ab}+O(\calh^+)),\bar{\calm}_{ab}^u =
\calm_{ab}^u\oplus\C(\Lambda_{ab}+O(\calh^+))$ and $\bar{\calm}_{ab}^t =
\calm_{ab}^t\oplus\C(\Lambda_{ab}+O(\calh^+))$.

%Let $\calp$ be the set of an ordered pair $\alpha = \{a,b\}$ such that
%$a\neq b$. We define $E_\alpha^u = E_{ab}^u, E_\alpha^t = E_{ab}^t$ and
%$\Lambda_\alpha = \Lambda_{ab}$. In the following, we consider the product
%$E_\alpha^uE_\beta^u$ etc. for $\alpha,\beta\in\calp$. The case decomposes
%into either $\alpha\cap\beta=\emptyset, \alpha = \beta$ or $\sharp
%(\alpha\cap\beta) =1$.

Recall that $[u]=u+O(\calh^+)$ is the image of $u\in \calh^+
$ in $A(\calh^+).$
\begin{lemma}\label{lemma:5.3.1a}
For distinct $a,b,c$ and $d$, 

\noindent
{\rm (i)} 
$\calm_{ab}*\bar{\calm}_{cd} = 0$, $\bar{\calm}_{cd}*\calm_{ab} =
0$.

\noindent
{\rm (ii)}
$[\Lambda_{ab}]*[\Lambda_{cd}] = [S_{ab}]*[S_{cd}].$
\end{lemma}
\begin{proof}
Let $u_{ab}\in\bar{\calm}_{ab}$ and $u_{cd}\in\bar{\calm}_{cd}$.  Remark
\ref{remark:4.1.9} shows that there exists a scalar $x$ such that
\begin{equation}\label{eqn:5.3.1a}
[u_{ab}]*[u_{cd}] = x[S_{ab}][S_{cd}].
\end{equation}
If $u_{ab}\in \calm_{ab}$, the left hand side of (\ref{eqn:5.3.1a}) acts
trivially on the top level of the module $M(1,\lambda)$. So the evaluation
of (\ref{eqn:5.3.1a}) on this top level shows $x=0$, which proves
$\calm_{ab}*\bar{\calm}_{cd} = 0$. Similarly,
$\bar{\calm}_{cd}*\calm_{ab} = 0$.

Next, the evaluation of $[\Lambda_{ab}]*[\Lambda_{cd}] = x[S_{ab}]*[S_{cd}]$ 
on the top level of the module $M(1,\lambda)$
leads $x=1$ as $\Lambda_{ab} =S_{ab}= \lambda_a\lambda_b.$
\end{proof}

As one can guess from Table 1, 2 and 3, the element $\Lambda_{ab}$
commutes with $\omega_a$ and $J_a$.

\begin{lemma}\label{lemma:5.3.2aa}
For distinct $a,b$ and any $c$,
$[\Lambda_{ab},\omega_c] =[\Lambda_{ab}, J_c] = 0$
in $A(\calh^+).$
\end{lemma}
\begin{proof}
In the case $a\neq c$, $b\neq c$, it is clear. Suppose $a=b$ or
$a=c$, then since $[\omega_c, \Lambda_{ab}] \sim 
(L_c(-1)+L_c(0))\Lambda_{ab}$ 
by Proposition \ref{proposition:3.2.1} (iii), the commutator belongs to
$\bar{\calm}_{ab}$. Since $[\omega_c, \Lambda_{ab}]$ acts trivially
on the top levels of the known irreducible modules by Tables 2 and
3, Lemma \ref{lemma:5.2.1a} asserts that
$[\Lambda_{ab},\omega_c] =0.$ The relation $[\Lambda_{ab},J_c] = 0$ is proved
by the same argument in the proof of Proposition
\ref{proposition:5.2.3}. That is, we  express $[\Lambda_{ab},J_c]$
as a linear combination of $E_{ab}^u, E_{ba}^u, E_{ab}^t, E_{ba}^t$
and $\Lambda_{ab}*\omega_c^j$ and then evaluate the expression on
the top levels  of the known irreducible modules to show that
$[\Lambda_{ab},J_c] = 0.$ 
\end{proof}

The following lemma is straightforward.

\begin{lemma}\label{lemma:5.3.1}
For all distinct $a,b$ and $c$,
\[
\begin{split}
S_{ab}(1,m)*S_{ac}(1,n) &= 2\omega_a*S_{bc}(m,n)
+\frac{1}{2}m(m+1)S_{bc}(m+2,n)\\ 
&\quad +m(m+1)S_{bc}(m+1,n)
+  \frac{1}{2}m(m+1)S_{bc}(m,n).
\end{split}
\]
\end{lemma}

\begin{lemma} \label{lemma:5.3.3a} Let $\{a,b\}\cap \{c,d\}$ have cardinality
one. Then we have  the following in $A(\calh^+).$

\noindent
{\rm (i)}
$\calm_{ab}^u*\calm_{cd}^t = \calm_{cd}^t *\calm_{ab}^u = 0$ and
$\Lambda_{ab}*\calm_{cd} = \calm_{cd}*\Lambda_{ab} = 0$.

\noindent
{\rm (ii)}
$E_{ab}^u*E_{cd}^u = \delta_{bc}E_{ad}^u$ and $E_{ab}^t*E_{cd}^t =
\delta_{bc}E_{ad}^t$.

Further suppose $a,b$ and $c$ are distinct, then

\noindent
{\rm (iii)}
$\Lambda_{ab}*\Lambda_{bc} = 2\omega_b*\Lambda_{ac}$.
\end{lemma}
\begin{proof}
We first consider part (i). Let $u_{ab}\in \bar{\calm}_{ab}$ and
$u_{cd}\in \bar{\calm}_{cd}$. We can assume $b=d$ without loss of
generality.  Then by Lemmas
\ref{lemma:4.2.3} and \ref{lemma:5.3.1} we have that
\[
%\label{1}
u_{ab}*u_{cd} = x_1\omega_b*S_{ac}(1,1) + x_2\Lambda_{ac} + y_1E_{ac}^u +
y_2E_{ca}^u + y_3E_{ac}^t + y_4E_{ca}^t
\]
with some scalars $x_i$ and  $y_j$. Note that  in part (i) case, 
element $u_{ab}*u_{cd}$ acts trivially on the top levels
of the modules $M(1,\lambda)$, $\calh^-$ and $\calh(\theta)^-.$
This implies that all $x_i$ and $y_j$ are zero.

Let us now prove the part (ii). Using the same argument in the
previous paragraph shows $E_{ab}^*E_{cd}^*\in \calm_{ad}$ if $b=c$ and 
$E_{ab}^*E_{cd}^*\in \calm_{bc}$ if $b\neq c$ as  $E_{ab}^u*E_{cd}^u$ acts
trivially on the 
top level of $M(1,\lambda)$. Then the result follows from Lemma
\ref{lemma:5.2.1a}.

Finally the evaluation of the equality
\[
\Lambda_{ab}*\Lambda_{bc} = x_1\omega_b*S_{ac}(1,1) + x_2\Lambda_{ac} +
y_1E_{ac}^u +
y_2E_{ca}^u + y_3E_{ac}^t + y_4E_{ca}^t
\]
on the top level of the modules $\calh^-, \calh(\theta)^-$ and $M(1,\lambda)$
(see Tables 1, 2, 3)  gives the linear system
\begin{align*}
&y_1= y_2= 0,\\
&x_1+32y_3 = 0, x_1+32y_4 = 0,\\
&x_1 = 2, x_2 = 0.
\end{align*}
Thus we have proved 
\begin{equation}\label{eqn:5.3.3a}
\Lambda_{ab}*\Lambda_{bc} = 2\omega_b*S_{ac}(1,1) -\frac{1}{16}E_{ac}^t
-\frac{1}{16}E_{ca}^t.
\end{equation}
Then substituting $S_{ac}(1,1) = E_{ac}^u+E_{ca}^u +\Lambda_{ac} +
\frac{1}{2}E_{ac}^t +\frac{1}{2}E_{ca}^t$ (see (\ref{eqn:n24}))
into (\ref{eqn:5.3.3a}) and using  Lemma \ref{lemma:5.2.2a}
gives  (iii).
\end{proof}

\begin{corollary}\label{corollary:5.3.5a}
For all distinct $a,b$ and $c$, $J_a*E_{bc}^u = E_{bc}^u*J_a = 0$ and
$J_a*E_{bc}^t = E_{bc}^t*J_a = \frac{1}{128}E_{bc}^t$ hold in
$A(\calh^+).$ 
\end{corollary}
\begin{proof}
By Lemma \ref{lemma:5.3.3a} (ii), $J_a*E_{bc}^* = J_a*E_{ba}^**E_{ac}^*.$
We have already proved  in Proposition \ref{proposition:5.2.3} that
 $J_a*E_{ba}^*\in \calm_{ba}^*.$ Using Lemma \ref{lemma:5.3.3a} (ii)
again shows that $J_a*E_{ba}^**E_{ac}^*\in \calm_{bc}^*.$ Similarly, 
$E_{bc}^**J_a =E_{ba}^**E_{ac}^**J_a\in\calm_{bc}^*$. Since
$J_a*E_{bc}^u,E_{bc}^u*J_a, J_a*E_{bc}^t-\frac{1}{128}E_{bc}^t$ and
$E_{bc}^t*J_a-\frac{1}{128}E_{bc}^t $ act trivially on the top levels of
$\calh^-$ and
$\calh(\theta)^-$, the corollary  follows. 
\end{proof}

We remark that the Corollary \ref{corollary:5.3.5a} completes the proof
of Proposition \ref{proposition:5.2.3}.

\begin{lemma}\label{lemma:5.3.4a}
For distinct $a$ and $b$, 
$(E_{ab}^u)^2 = (E_{ab}^t)^2 = 0, E_{ab}^uE_{ab}^t = E_{ab}^tE_{ab}^u
= 0$ and $E_{ab}^u\Lambda_{ab} = \Lambda_{ab}E_{ab}^u = 0,
E_{ab}^t\Lambda_{ab} = \Lambda_{ab}E_{ab}^t = 0$ hold in $A(\calh^+).$  
\end{lemma}
\begin{proof}
Since $S_{ab}(m,n)*S_{ab}(s,t)\in \calw$ for all $m,n,s,t>0$ we see 
that $E_{ab}^*E_{ab}^*$ is a polynomial in 
$\omega_a,\omega_b,J_a$ and $J_b$ (cf.~Corollary
\ref{corollary:4.1.2}).  In particular, $E_{ab}^*E_{ab}^*$ commutes with
$\omega_a$. Then by Lemma \ref{lemma:5.2.2a}, we have 
\[
0=[\omega_a,E_{ab}^uE_{ab}^u] =\omega_a
E_{ab}^uE_{ab}^u-E_{ab}^uE_{ab}^u\omega_a
=E_{ab}^uE_{ab}^u.
\]
Similarly, $(E_{ab}^t)^2=E_{ab}^uE_{ab}^t = E_{ab}^tE_{ab}^u =
0$.

Since $\Lambda_{ab}$ commutes with $\omega_a$ by Lemma \ref{lemma:5.3.2aa},
we also have $E_{ab}^u\Lambda_{ab}\omega_a = 0$. Thus the same argument in
the previous paragraph shows that  $E_{ab}^u\Lambda_{ab}= 0$. The rest of
the equalities are proved
similarly.
\end{proof}

The following is an immediate corollary of Lemma~\ref{lemma:5.3.3a} (ii).

\begin{lemma}\label{lemma:5.3.5a}
For distinct $a,b$ and $c$, $E_{ab}^uE_{ba}^u = E_{ac}^uE_{ca}^u$ and
$E_{ab}^tE_{ba}^t = E_{ac}^tE_{ca}^t$.
\end{lemma}

Thanks to Lemma \ref{lemma:5.3.4a}, we define new elements $E_{aa}^* =
E_{ab}^*E_{ba}^*$ for all $1\leq a\leq \ell$, which is independent of the
choice of the index $b$. 

\begin{lemma}\label{lemma:5.3.6a}
For any $a$ and distinct $b,c$, 

\noindent
{\rm (i)}
$E_{aa}^uE_{bc}^u = \delta_{ab}E_{ac}^u$ and $E_{bc}^uE_{aa}^u =
\delta_{ca}E_{ba}^u$.

\noindent
{\rm (ii)}
$E_{aa}^tE_{bc}^t = \delta_{ab}E_{ac}^t$ and $E_{bc}^tE_{aa}^t =
\delta_{ca}E_{ba}^t$.

\noindent
{\rm (iii)}
$E_{aa}^uE_{bc}^t = E_{bc}^tE_{aa}^u = 0$.

\noindent
{\rm (iv)}
$E_{aa}^tE_{bc}^u = E_{bc}^uE_{aa}^t = 0$.

\noindent
{\rm (v)}
$E_{aa}^u\Lambda_{bc} =\Lambda_{bc}E_{aa}^u = 0$.

\noindent
{\rm (vi)}
$E_{aa}^t\Lambda_{bc} =\Lambda_{bc}E_{aa}^t = 0$.
\end{lemma}
\begin{proof}
Suppose $a\neq b$. Then Lemma \ref{lemma:5.3.3a} (ii) or Lemma
\ref{lemma:5.3.4a} shows $E_{aa}^*E_{bc}^* = E_{ab}^*E_{ba}^*E_{bc}^* = 0$.

If $a= b$, then $E_{aa}^*E_{bc}^* = E_{ac}^*E_{ca}^*E_{ac}^*$. Since
$E_{ac}^*E_{ca}^*$ lies in the subalgebra generated by 
$\omega_a,\omega_c, J_a,J_c$, we use Lemma
\ref{lemma:5.2.2a} and Proposition \ref{proposition:5.2.3} (also
see  Corollary \ref{corollary:5.3.5a}) to show $E_{aa}^*E_{ac}^*\in
\calm_{ac}^*$. Then the evaluation method proves $E_{aa}^*E_{ac}^* =
E_{ac}^*$.

The rest of the equalities in (i), (ii) and (iii)--(vi) are proved
similarly.
\end{proof}

Now we go back to $E_{ab}^uE_{ba}^t$ and $E_{ab}^tE_{ba}^u$.

\begin{lemma}\label{lemma:5.3.7a}
Suppose $a\neq b$. Then $E_{ab}^uE_{ba}^t = E_{ab}^tE_{ba}^u
= 0$.
\end{lemma}
\begin{proof}
Since $E_{ab}^* = E_{ab}^*E_{bb}^*$ by Lemma \ref{lemma:5.3.6a} (i) and
(ii), we have $E_{ab}^uE_{ba}^t = E_{ab}^uE_{bb}^uE_{ba}^t = 0$ and
$E_{ab}^tE_{ba}^u = E_{ab}^tE_{bb}^tE_{ba}^u = 0$
by Lemma \ref{lemma:5.3.6a} (iii) and (iv).
\end{proof}

Using the same argument in the proof of Lemma \ref{lemma:5.3.7a}, 
we prove:

\begin{lemma}\label{lemma:5.3.10a}
For any distinct $a$ and $b$, $E_{aa}^uE_{bb}^u =
E_{aa}^tE_{bb}^t = 0$ and
$E_{aa}^uE_{bb}^t = E_{bb}^tE_{aa}^u = 0$.
\end{lemma}

\begin{lemma}\label{lemma:5.3.11a}
For any $a$, $(E_{aa}^u)^2 = E_{aa}^u, (E_{aa}^t)^2 = E_{aa}^t$ and
$E_{aa}^uE_{aa}^t = E_{aa}^tE_{aa}^u = 0$.
\end{lemma}
\begin{proof}
Using Lemma \ref{lemma:5.3.6a} (i) and (ii), we see $E_{aa}^*E_{aa}^* =
E_{ab}^*E_{ba}^*E_{aa}^*= E_{ab}^*E_{ba}^* = E_{aa}^*$.
Similarly we prove $E_{aa}^uE_{aa}^t =E_{aa}^tE_{aa}^u = 0$ by Lemma
\ref{lemma:5.3.6a} (iii) and (iv).
\end{proof}

Summarizing  we have: 

\begin{proposition}\label{proposition:5.3.3}
For any $a,b,c$ and $d$, 

\noindent
{\rm (i)} 
$E_{ab}^uE_{cd}^u = \delta_{bc}E_{ad}^u$.

\noindent
{\rm (ii)}
$E_{ab}^tE_{cd}^t = \delta_{bc}E_{ad}^t$.

\noindent
{\rm (iii)}
$E_{ab}^uE_{cd}^t = 0,E_{cd}^tE_{ab}^u=0$.

\noindent
{\rm (iv)}
$\Lambda_{ab}E_{cd}^u = E_{cd}^u\Lambda_{ab} = 0$ and
$\Lambda_{ab}E_{cd}^t = E_{cd}^t \Lambda_{ab}= 0,\,(a\neq b)$.

Further suppose $a,b,c$ distinct,

\noindent
{\rm (v)}
$\Lambda_{ab}\Lambda_{bc}= 2\omega_b*\Lambda_{ac}.$

\end{proposition}

The following corollary follows from Remark \ref{remark:5.1.3},
Lemma \ref{lemma:5.3.2aa} and Proposition
\ref{proposition:5.3.3}.

\begin{corollary}\label{c1}
For distinct $a$ and $b$, $\Lambda_{ab}$ lies in the the center of
$A(\calh^+)$.
\end{corollary}

Let $\cali^u =
\bigl(\oplus_{\substack{1\leq a,
b\leq \ell,\\
a\neq b}}\calm_{ab}^u\bigr)\bigoplus\bigl(\oplus_{1\leq a\leq \ell}\C
E_{aa}^u\bigr),\cali^t  =
\bigl(\oplus_{\substack{1\leq a,
b\leq \ell,\\
a\neq b}}\calm_{ab}^t\bigr)\bigoplus\bigl(\oplus_{1\leq a\leq \ell}\C
E_{aa}^t\bigr)$ and $\cali =
\cali^u\oplus\cali^t$.

\begin{proposition}\label{proposition:5.3.4}
The vector spaces $\cali^u$ and $\cali^t$ are ideals of
$A(\calh^+)$ and are isomorphic to the matrix algebra
$\operatorname{Mat}_\ell(\C)$ as associative algebras. Furthermore,
$\cali^u*\cali^t =\cali^t* \cali^u=0$. 
\end{proposition}
\begin{proof}
It suffice to prove that $\cali^u$ and $\cali^t$ are ideals. But this is
clear from Remark \ref{remark:5.1.3}, 
Lemma \ref{lemma:5.2.2a}, Propositions \ref{proposition:5.2.3}
and \ref{proposition:5.3.3}.
\end{proof}

By Corollary \ref{corollary:4.1.2}, $\omega_a,\omega_b,J_a, J_b$ 
commute in $A(\calh^+).$ Combining Remark \ref{remark:5.1.3},
Corollary \ref{c1} and Proposition \ref{proposition:5.3.4}
together gives the following result. 

\begin{proposition}\label{proposition:5.3.5}
The quotient algebra $A(\calh^+)/\cali$ is abelian and is generated by the
image of $\omega_a,J_a$ and $\Lambda_{ab}$ for $1\leq a,b\leq \ell$.
\end{proposition}

\section{Classification of irreducible modules for $\calh^+$}
\label{bsection:6}

In this conclusion section we first give more relations among generators
of $A(\calh^+)$ and classify the irreducible modules 
for $\calh^+$ by classifying the simple modules for $A(\calh^+).$

\subsection{Relations among $\omega_a, J_a$ and $E_{ab}^*$}
\label{subsection:6.1}

Set $H_a = J_a + \omega_a -4\omega_a^2,\,(i\leq a\leq\ell)$. Let us
recall that
\[
S_{ab}(m,n) = h_a(-m)h_b(-n)\unit
\]
for $m,n>0.$ In this subsection we shall compute
$[S_{ab}(1,m),S_{ab}(1,n)]$ for certain $m,n>$ in two ways and obtain
new relations in $A(\calh^+)$. Since both $S_{ab}(1,m)$ and
$S_{ab}(1,n)$ are linear combinations of $E_{ab}^u, E_{ba}^u, E_{ab}^t,
E_{ba}^t$ and $\Lambda_{ab}$, $[S_{ab}(1,m),S_{ab}(1,n)]$ is a linear
combinations of $E_{aa}^u, E_{bb}^u, E_{aa}^t, E_{bb}^t$ on the one
hand. On the other hand,  a direct computation of
$[S_{ab}(1,m),S_{ab}(1,n)]$ by definition shows that
$[S_{ab}(1,m),S_{ab}(1,n)]$ is a polynomial in 
$\omega_a,\omega_b,J_a,J_b.$ This will give two new relations. The
other two new relations are derived in  different ways.

Recall from Remark \ref{remark:5.1.2} that
\begin{align*}
S_{ab}(1,1)& = E_{ab}^u + E_{ba}^u +\Lambda_{ab}+\frac{1}{2}E_{ab}^t +
\frac{1}{2}E_{ba}^t,\\
S_{ab}(1,2) &= -2E_{ab}^u- \Lambda_{ab}-\frac{3}{4}E_{ab}^t -
\frac{1}{4}E_{ba}^t,\quad\text{and}\\
S_{ab}(1,4)&= -4 E_{ab}^u- \Lambda_{ab} - \frac{35}{32}E_{ab}^t -
\frac{5}{32}E_{ba}^t.
\end{align*}

\begin{lemma}\label{l6.1} If $a\neq b$, then 
\begin{align*}
&-\frac{2}{9}H_a + \frac{2}{9}H_b =2E_{aa}^u -2E_{bb}^u
+\frac{1}{4}E_{aa}^t -\frac{1}{4}E_{bb}^t,\quad\text{and}\\
&-\frac{4}{135}(2\omega_a + 13)H_a + \frac{4}{135}(2\omega_b + 13) H_b 
 =4E_{aa}^u -4E_{bb}^u +\frac{15}{32}E_{aa}^t
-\frac{15}{32}E_{bb}^t.
\end{align*}
\end{lemma}
\begin{proof} First, we see from Proposition \ref{proposition:3.1.2} (iii) 
that 
\[[S_{ab}(1,1),S_{ab}(1,m)]\sim
S_{ab}(1,1)_{0}S_{ab}(1,m)+S_{ab}(1,1)_{1}S_{ab}( 1,m).
\]
 Note that
\[
S_{ab}(1,1)_{k} = \sum_{i+j=k-1}\NO h_a(i)h_b(j)\NO
\]
and $i$ or $j$ is nonnegative if $k\geq 0$. Thus
$S_{ab}(1,1)_{0}S_{ab}(1,m)+S_{ab}(1,1)_{1}S_{ab}( 1,m)\in \calh^+_{a}
+\calh^+_{b}.$ 

We now deal with $[S_{ab}(1,1),S_{ab}(1,2)].$
 From Proposition \ref{proposition:5.3.3} and Corollary
\ref{c1}, 
\[
[S_{ab}(1,1),S_{ab}(1,2)]= 2E_{aa}^u -2E_{bb}^u +\frac{1}{4}E_{aa}^t
-\frac{1}{4}E_{bb}^t.
\]
Since the weights of 
$S_{ab}(1,1)_{0}S_{ab}(1,2)$ and $S_{ab}(1,1)_{1}S_{ab}( 1,2)$ are less
than or equal to 4, we 
see that \begin{equation}\label{eqn:39}
[S_{ab}(1,1),S_{ab}(1,2)] = x H_a + f(\omega_a) + y H_b +
g(\omega_b)
\end{equation}
with some scalars $x, y$ and polynomials $f$ and $g.$

Evaluating (\ref{eqn:39}) on the top level of the module $M(1,\lambda)$
and  noting that $H_a=H_b=0$ (see Table 3) yields
$f(\frac{1}{2}\lambda_a^2) + g(\frac{1}{2}\lambda_b^2) = 0$
for all $0\neq \lambda \in \frakh$, which implies $f= g = 0$. Therefore, we
have
\begin{equation}\label{eqn:40}
x H_a + y H_b =2E_{aa}^u -2E_{bb}^u
+\frac{1}{4}E_{aa}^t-\frac{1}{4}E_{bb}^t.
\end{equation}
Further the evaluation of (\ref{eqn:40}) on the top level of the
module $\calh(\theta)^-$ shows that $x = -2/9$ and $y = 2/9.$

Similarly, compute $[S_{ab}(1,1),S_{ab}(1,4)]$ in two different
ways gives 
\[
 (x\omega_a+y)H_a+f(\omega_a) + (z\omega_b+w)H_b
+f(\omega_b)=4E_{aa}^u -4E_{bb}^u +\frac{15}{32}E_{aa}^t
-\frac{15}{32}E_{bb}^t
\]
for some scalars $x,y,z,w$ and polynomials $f$ and $g.$ 
Again the evaluation on
$M(1,\lambda)$ shows that $f=g=0.$ That is 
\begin{equation}\label{eqn:41}
 (x\omega_a+y)H_a + (z\omega_b+w)H_b =4E_{aa}^u -4E_{bb}^u 
+\frac{15}{32}E_{aa}^t
-\frac{15}{32}E_{bb}^t.
\end{equation}

Now the evaluation of (\ref{eqn:41}) on the top level of the module $\calh^-$
shows
\begin{equation}\label{eqn:42}
x+y = -\frac{4}{9}\quad \text{and}\quad z+w = \frac{4}{9},
\end{equation}
%and the evaluation of (\ref{eqn:41}) on the top level of the module
%$\calh(\theta)^+$ gives
%\begin{equation}\label{eqn:43}
%x+16y + z+16w = 0.
%\end{equation}
 the evaluation of (\ref{eqn:41}) on the top level of the module
$\calh(\theta)^-$ yields
\begin{equation}\label{eqn:44}
-\frac{153}{256}x-\frac{9}{8}y = \frac{15}{32}\quad\text{and}\quad
-\frac{153}{256}z-\frac{9}{8}w = -\frac{15}{32}.
\end{equation}
Solving the linear system (\ref{eqn:42}) and (\ref{eqn:44})
we find
\[
x = -\frac{8}{135},\quad y = -\frac{52}{135},\quad z =
\frac{8}{135},\quad w = \frac{52}{135}.
\]
\end{proof}

We need two more relations.
\begin{lemma}\label{ll1}
If $a\ne b$, then 
\begin{align*}
&\omega_b H_a = -\frac{2}{15}(\omega_a -1)H_a + \frac{1}{15}(\omega_b
-1)H_b,\quad\text{and}\\
&\Lambda_{ab}^2 =
4\omega_a\omega_b-\frac{1}{9}(H_a+H_b) - (E_{aa}^u+E_{bb}^u)
-\frac{1}{4}(E_{aa}^t+E_{bb}^t).
\end{align*}
\end{lemma}
\begin{proof}
Recall that 
\begin{align*}
\omega_a &= \frac{1}{2}h_a(-1)^2,\\
\omega_a^2 &\sim (L_a(-2)+L_a(-1))\omega_a= \frac{1}{4}h_a(-1)^4\unit +
h_a(-3)h_a(-1)\unit + h_a(-2)h_a(-1)\unit,\\ 
J_a &= h_a(-1)^4\unit - 2h_a(-3)h_a(-1)\unit +
\frac{3}{2}h_a(-2)^2\unit.
\end{align*}
So 
\[
H_a \sim -6h_a(-3)h_a(-1)\unit - 4h_a(-2)h_a(-1)\unit +
\frac{3}{2}h_a(-2)^2\unit+ \frac{1}{2}h_a(-1)^2\unit.
\]
Now using $L(-1)S_{aa}(2,1)+ L(0)S_{aa}(2,1)\sim 0$
(cf. Proposition \ref{proposition:3.1.2} (iv)),
we have 
\[
2h_a(-3)h_a(-1)\unit + h_a(-2)^2\unit + 3h_a(-2)h_a(-1)\unit\sim 0.
\]
and
\[
H_a \sim -9S_{aa}(1,3) -\frac{17}{2}S_{aa}(1,2) +\frac{1}{2}S_{aa}(1,1).
\]
Then Lemma \ref{lemma:4.2.4} shows
\begin{equation}\label{eqn:45}
\omega_bH_a = c\omega_a\omega_b + (\alpha \omega_a +\beta)H_a + (x \omega_b
+ y) H_b +f(\omega_a) + g(\omega_b)
\end{equation}
with some scalars $c,\alpha,\beta,x,y$ and polynomials $f$ and $g$.

Now the evaluation of (\ref{eqn:45}) on the top level of the module
$M(1,\lambda)$ shows
\[
\frac{c}{4}\lambda_a^2\lambda_b^2 + f(\frac{1}{2}\lambda_a^2) +
g(\frac{1}{2}\lambda_b^2) = 0.
\]
This implies $c = 0$ and $f=g = 0$.
Thus we have 
\begin{equation}\label{eqn:46}
\omega_b H_a = (\alpha\omega_a+\beta)H_a + (x\omega_b +y)H_b.
\end{equation}
Finally the evaluation of (\ref{eqn:46}) on the top levels of the modules
$\calh^-$ and $\calh(\theta)^-$ respectively gives
\begin{align}
&\alpha +\beta = 0,\quad x+y = 0,\label{eqn:47}\\
%&\alpha + 16\beta + x + 16y = 1,\label{eqn:48}\\
&-\frac{153}{256}\alpha -\frac{9}{8}\beta = -\frac{9}{128},\quad
-\frac{153}{256}x -\frac{9}{8}y = \frac{9}{256}.\label{eqn:49}
\end{align}
Then solving the linear system (\ref{eqn:47}) and (\ref{eqn:49}) gives
 $\alpha = -2/15, \beta = 2/15, x = 1/15$ and $y =
-1/15$.

For the second relation, note from Lemma \ref{lemma:5.1.1}
and Proposition \ref{proposition:5.3.3}
that $\Lambda_{ab} = \Lambda_{ba}$ and $\Lambda_{ab}E_{cd}^u =
\Lambda_{ab}E_{cd}^t= 0.$ Thus
\[
S_{ab}(1,1)^2 = E_{aa}^u + E_{bb}^u + \Lambda_{ab}^2 +
\frac{1}{4}E_{aa}^t +
\frac{1}{4}E_{bb}^t.
\]
On the other hand, the evaluation method shows $S_{ab}(1,1)^2 =
4\omega_a\omega_b - \frac{1}{9}H_a -\frac{1}{9}H_b$
and then we have
\[
\Lambda_{ab}^2 =
4\omega_a\omega_b-\frac{1}{9}(H_a+H_b) - (E_{aa}^u+E_{bb}^u)
-\frac{1}{4}(E_{aa}^t+E_{bb}^t).
\]
\end{proof}

Summarizing, we have: 

\begin{proposition}\label{proposition:6.2.1}
For distinct $a$ and $b$,
\begin{align}
&(70H_a+1188\omega_a^2-585\omega_a +27)H_a = 0,\label{eqn:50}\\
&(\omega_a-1) \bigl(\omega_a-\frac{1}{16}\bigr)
\bigl(\omega_a-\frac{9}{16}\bigr) H_a= 0,\label{eqn:51}\\
&-\frac{2}{9}H_a + \frac{2}{9}H_b =2E_{aa}^u -2E_{bb}^u +\frac{1}{4}E_{aa}^t
-\frac{1}{4}E_{bb}^t,\label{eqn:52}\\
&-\frac{4}{135}(2\omega_a + 13)H_a + \frac{4}{135}(2\omega_b + 13) H_b =
4E_{aa}^u -4E_{bb}^u +\frac{15}{32}E_{aa}^t
-\frac{15}{32}E_{bb}^t,\label{eqn:53}\\
&\omega_b H_a = -\frac{2}{15}(\omega_a -1)H_a + \frac{1}{15}(\omega_b
-1)H_b,\label{eqn:54}\\ 
&\Lambda_{ab}^2 =
4\omega_a\omega_b-\frac{1}{9}(H_a+H_b) - (E_{aa}^u+E_{bb}^u)
-\frac{1}{4}(E_{aa}^t+E_{bb}^t),\label{eqn:55}\\ 
&\Lambda_{ab}\Lambda_{bc}= 2\omega_b*\Lambda_{ac}\quad
\text{for $a,b,c$ distinct}.\label{eqn:56}
\end{align}
\end{proposition}

(\ref{eqn:50}) and (\ref{eqn:51}) were given in \cite{DN1} by noting
that $70H_a+1188\omega_a^2-585\omega_a +27=70J_a+908\omega^2-515\omega+27.$

\subsection{Irreducible modules for $A(\calh^+)$ and the classification
result }
\label{subsection:6.2}

Thanks to the relations we have obtained in the preceding subsections we are
able to classify all the irreducible modules for the algebra
$A(\calh^+)$ and then for the vertex operator algebra $\calh^+$.

Recall $\cali = \cali^u\oplus\cali^t$. Let $W$ be an irreducible 
$A(\calh^+)$-module. There are two cases: $\cali W\ne 0$ or 
$\cali W=0.$ If $\cali W\ne 0$ then either $\cali^u W \ne 0$ 
or $\cali^t W\ne 0.$ If  $\cali^u W\ne 0$ then we must have 
$\cali^u W=W$ and $\cali^t W=0$ as $\cali^u, \cali^t $ are
ideals of $A(\calh^+)$ and  $\cali^u\cap\cali^t=0.$ Thus
$W$ is a simple module for $\cali^u$ and is isomorphic to
the top level of $\calh^-$ as $A(\calh^+)$-modules. Similarly,
if $\cali^t W\ne 0,$ then $W$ is isomorphic to the top level
of $\calh(\theta)^-.$

Now we assume that $\cali W=0.$ Then $W$ is a module for 
quotient algebra $A(\calh^+)/\cali$  which 
is commutative by Proposition \ref{proposition:5.3.5}.
Since $A(\calh^+)$ has countable dimension, $W$ is $1$-dimensional.
Then each element in $A(\calh^+)$ acts a scalar on $W.$

Suppose $H_a= 0$ for some $a$ on $W.$ Then (\ref{eqn:52}) shows $H_a = 0$ for
all $a$. Let $\omega_a = \frac{1}{2}\lambda_a^2$ with $\lambda_a\in\C$.
Then we have $J_a = \lambda_a^4 -\frac{1}{2}\lambda_a^2$ as $H_a = 0$.
Now (\ref{eqn:55}) shows
\[
\Lambda_{ab}^2 = 4\omega_a\omega_b = \lambda_a^2\lambda_b^2.
\]
on $W.$ 
Therefore $\Lambda_{ab} = \varepsilon_{ab}\lambda_a\lambda_b$ where
$\varepsilon_{ab} = \pm 1$. Substituting the relation into (\ref{eqn:56}), 
we have $\varepsilon_{ab}\varepsilon_{bc} =
\varepsilon_{ac}$. Therefore there exists a map $\varepsilon$
from the index set $\{1,2,\dots, \ell\,\}$ to $\{\pm 1\}$ such that
$\varepsilon_{ab} = (-1)^{\varepsilon(a)-\varepsilon(b)}$. Then replacing
$\lambda_a$ by $(-1)^{\varepsilon(a)}\lambda_a$, we have
\[
\omega_a  = \frac{1}{2}\lambda_a^2,\quad J_a = \lambda_a^4
-\frac{1}{2}\lambda_a^2,\quad \text{and}\quad \Lambda_{ab} =
\lambda_a\lambda_b.
\]
Thus the module $W$ is isomorphic to the top level of the module
$M(1,\lambda)$ as $A(\calh^+)/\cali$ is generated by $\omega_a, J_a$ and
$\Lambda_{ab}$ (see  Proposition \ref{proposition:5.3.5}.)

Next suppose $H_a\neq 0$ on $W$ for some $a$. Then (\ref{eqn:52}) says $H_a =
H_b  $ for all $a$ and $b$. Then (\ref{eqn:53}) shows
$\omega_a-\omega_b = 0$ on $W.$ Now set $\omega_a =\lambda$. Then 
(\ref{eqn:54}) implies $\lambda = 1/16$ and then $J_a = 3/128$ 
by (\ref{eqn:50}). Finally (\ref{eqn:55}) shows
$\Lambda_{ab}^2 = 0$, i.e., $\Lambda_{ab} = 0$. Thus the
module $W$ in this case 
is isomorphic to the top level of the module $\calh(\theta)^+$.

Note that the top level of $M(1,0)$ is also the top level of
$\calh^+.$ 
Thus we have proved:
\begin{proposition}\label{proposition:7.1.1}
Any irreducible module for $A(\calh^+)$ is isomorphic
either to the top level of the module $M(1,\lambda),\,(0\ne \lambda\in \frakh)$
or to the top level of $\calh^+$ or to the top level of $\calh(\theta)^+$.
\end{proposition}

Finally by Theorem \ref{theorem:3.1.1} (iii) we have
proved:

\begin{theorem}\label{theorem:7.2.2}
Any irreducible admissible 
module for the vertex operator algebra $\calh^+$ is
isomorphic to one of the following modules;
\[
\calh^+,\quad \calh^-,\quad M(1,\lambda)\simeq M(1,-\lambda),\,(\lambda\neq
0\in \frakh),\quad
\calh(\theta)^+,\quad \calh(\theta)^-.
\]
In particular, any irreducible admissible module for $\calh^+$ is an
ordinary module.
\end{theorem}

%References

\end{document}